\documentclass[final,1p,times]{elsarticle}

\usepackage{amssymb}
 \usepackage{amsthm}
\usepackage{amscd}
\usepackage{amsmath}
\usepackage{amsfonts}
\usepackage{amssymb}
\usepackage{graphicx}
\usepackage{color}
\newtheorem{theorem}{Theorem}
\newtheorem{proposition}[theorem]{Proposition}
\newtheorem{lemma}[theorem]{Lemma}

\usepackage{mathrsfs}
\usepackage{titletoc}

\newcommand{\D}{\Delta}
\newcommand{\ra}{\rightarrow}
\newcommand{\p}{\partial}
\newcommand{\f}{\frac}

\newcommand{\be}{\begin{equation}}
\renewcommand{\ra}{\rightarrow}
\newcommand{\ee}{\end{equation}}
\newcommand{\bea}{\begin{eqnarray}}
\newcommand{\eea}{\end{eqnarray}}
\newcommand{\bna}{\begin{eqnarray*}}
\newcommand{\ena}{\end{eqnarray*}}

\renewcommand{\le}{\left}
\newcommand{\ri}{\right}

\journal{Ann. Sc. Norm. Sup. Pisa Cl. Sci.}

\begin{document}

\begin{frontmatter}

\title{Trudinger-Moser inequalities on a closed Riemannian surface with the action of a finite
isometric group}

\author{Yu Fang}
 \ead{2011202402@ruc.edu.cn}
\author{Yunyan Yang\footnote{Corresponding author}}
 \ead{yunyanyang@ruc.edu.cn}
 \address{Department of Mathematics,
Renmin University of China, Beijing 100872, P. R. China}

\begin{abstract}
Let $(\Sigma,g)$ be a closed Riemannian surface, $W^{1,2}(\Sigma,g)$ be the usual Sobolev space,
$\textbf{G}$ be a finite isometric group acting on $(\Sigma,g)$, and
$\mathscr{H}_\textbf{G}$ be the function space including all functions $u\in W^{1,2}(\Sigma,g)$ with
$\int_\Sigma udv_g=0$ and $u(\sigma(x))=u(x)$ for all $\sigma\in \textbf{G}$ and all $x\in\Sigma$. Denote the number of distinct
points of the set $\{\sigma(x): \sigma\in \textbf{G}\}$ by $I(x)$ and $\ell=\min_{x\in \Sigma}I(x)$.
Let $\lambda_1^\textbf{G}$ be the first eigenvalue of the Laplace-Beltrami operator on the space $\mathscr{H}_\textbf{G}$.
Using blow-up analysis,
we prove that if $\alpha<\lambda_1^\textbf{G}$ and $\beta\leq 4\pi\ell$, then there holds
$$\sup_{u\in\mathscr{H}_\textbf{G},\,\int_\Sigma|\nabla_gu|^2dv_g-\alpha \int_\Sigma u^2dv_g\leq 1}\int_\Sigma e^{\beta u^2}dv_g<\infty;$$
if $\alpha<\lambda_1^\textbf{G}$ and $\beta>4\pi\ell$, or  $\alpha\geq \lambda_1^\textbf{G}$ and $\beta>0$, then
the above supremum is infinity; if $\alpha<\lambda_1^\textbf{G}$ and $\beta\leq 4\pi\ell$, then the above supremum can be attained. Moreover,
similar inequalities involving higher order eigenvalues are obtained. Our results  partially
improve original inequalities of J. Moser \cite{Moser}, L. Fontana \cite{Fontana} and W. Chen \cite{Chen-90}.

\end{abstract}


\end{frontmatter}

\section{Introduction}
Let $\Omega\subset\mathbb{R}^n$ be a smooth bounded domain, $W_0^{1,n}(\Omega)$ be the usual Sobolev space,
and $\omega_{n-1}$ be the area of the unit sphere in $\mathbb{R}^n$. It was proved by Moser \cite{Moser}
that for any $\alpha\leq \alpha_n=n\omega_{n-1}^{1/(n-1)}$, there holds
\be\label{Moser-1}\sup_{u\in W_0^{1,n}(\Omega),\,\int_\Omega|\nabla u|^ndx\leq 1}\int_\Omega e^{\alpha |u|^{n/(n-1)}}dx<\infty.\ee
Moreover, $\alpha_n$ is the best constant in the sense that if $\alpha>\alpha_n$, the integrals in the above inequality are still finite,
but the supremum is infinity. Such kind of inequalities are known as Trudinger-Moser inequalities in literature. Earlier contributions
are due to Yudovich \cite{Yudovich}, Pohozaev \cite{Pohozaev}, Peetre \cite{Peetre} and Trudinger \cite{Trudinger}.
Let $\lambda_1(\Omega)$ be the first eigenvalue of the Laplace operator with respect to the Dirichlet boundary condition.
Adimurthi-Druet \cite{A-D}
proved that for any $\alpha<\lambda_1(\Omega)$, there holds
\be\label{A-D}\sup_{u\in W_0^{1,2}(\Omega),\,\int_\Omega|\nabla u|^2dx\leq 1}\int_\Omega e^{4\pi u^2(1+\alpha\|u\|_2^2)}dx<\infty;\ee
moreover, if $\alpha\geq \lambda_1(\Omega)$, then the above supremum is infinity, where $\|u\|_2^2=\int_\Omega u^2dx$.
The inequality (\ref{A-D}) is stronger than (\ref{Moser-1}) and was extended by the second named author \cite{Yang-JFA-06} to the
higher dimensional case.
Later, Tintarev \cite{Tintarev} proved among other results that for any $\alpha<\lambda_1(\mathbb{B}_R(0))$, there holds
\be\label{Tintarev}\sup_{u\in W_0^{1,2}(\Omega),\,\int_\Omega|\nabla u|^2dx-\alpha \int_\Omega u^2dx\leq 1}\int_\Omega e^{4\pi u^2}dx<\infty,\ee
where $\mathbb{B}_R(0)$ denotes the ball centered at $0$ with radius $R$ and its measure is equal to that of $\Omega$.
As one expected, $\lambda_1(\mathbb{B}_R(0))$ can be replaced by $\lambda_1(\Omega)$, which is a consequence of (\cite{Yang-JDE-15}, Theorem $1$).

One can ask whether the supremum in (\ref{Moser-1}) can be attained or not. Existence of extremal functions
was proved first by Carleson-Chang \cite{C-C}
in the case that $\Omega$ is the unit ball, then by Struwe \cite{Struwe} in the case that $\Omega$ is close to a ball in the sense of measure,
  later by  Flucher \cite{Flucher} when $\Omega$ is a planar domain, and finally by Lin
\cite{Lin} when $\Omega$ is a  domain in $\mathbb{R}^n$. In \cite{InterJM},
the second named author claimed that the supremum in (\ref{A-D}) can be attained for
all $0\leq\alpha<\lambda_1(\Omega)$. We remark that there is a mistake during that test function computation (\cite{InterJM}, page 338, line 8).
In fact, in two dimensions,  extremal function for (\ref{A-D}) exists only for sufficiently small $\alpha$, see for example \cite{Yang-Tran}.
 Concerning extremal functions for inequalities of the type (\ref{A-D}),
 we refer the reader to \cite{Lu-Yang-1,Lu-Yang-2,de-doO,YangJGA,Yang-Zhu-JFA,Iula-Mancini,ZhuJ,YuanHuang,YuanZhu,Nguyen}.
 As a comparison,  it was proved in \cite{Yang-JDE-15} that the supremum in (\ref{Tintarev}) can be attained for all $\alpha<\lambda_1(\Omega)$.
It is remarkable that (\ref{Tintarev}) is stronger than (\ref{A-D}), however, there is no relation on existence of
extremal functions between (\ref{A-D}) and (\ref{Tintarev}).

Let $(\mathbb{S}^2,g_0)$ be the 2-dimensional sphere $x_1^2+x_2^2+x_3^2=1$ with the metric $g_0=dx_1^2+dx_2^2+dx_3^2$
and the corresponding volume element $dv_{g_0}$. According to Moser \cite{Moser}, one can find a constant $C$ such that
for all functions $u$ with $\int_{\mathbb{S}^2}|\nabla_{g_0} u|^2dv_{g_0}\leq 1$ and $\int_{\mathbb{S}^2}udv_{g_0}=0$,
\be\label{Mosr-2}\int_{\mathbb{S}^2} e^{4\pi u^2}dv_{g_0}\leq C.\ee
Concerning all even functions $u$, it was indicated by Moser \cite{Moser-73} that the best constant $\alpha_2=4\pi$ would double.
Namely, there exists a constant $C$ such that for all functions $u$ satisfying $u(-x)=u(x)$, $\forall x\in \mathbb{S}^2$, $\int_{\mathbb{S}^2}|\nabla_{g_0}u|^2dv_{g_0}\leq 1$,
and $\int_{\mathbb{S}^2}udv_{g_0}=0$, there holds
\be\label{Mosr-even}
\int_{\mathbb{S}^2} e^{8\pi u^2}dv_{g_0}\leq C.\ee
Later, by using an isoperimetric inequality on closed Riemannian surfaces with conical singularities, Chen \cite{Chen-90}
proved a Trudinger-Moser inequality for a class of ``symmetric" functions, which particularly generalized
 (\ref{Mosr-2}) and (\ref{Mosr-even}).

Let $(M,g)$ be a closed $n$-dimensional Riemannian manifold. Among other results, it was proved by Fontana \cite{Fontana}
that there exists a constant
$C$, depending only on $(M,g)$, such that if $u\in W^{1,n}(M,g)$ satisfies $\int_M|\nabla_gu|^ndv_g\leq 1$ and $\int_Mudv_g=0$, then
\be\label{Fontana}\int_M e^{\alpha_n |u|^{n/(n-1)}}dv_g\leq C.\ee
The existence of extremal functions for (\ref{Fontana}) was obtained by Li \cite{Lijpde,Liscience}. Precisely, there exists some
$u_0\in W^{1,n}(M)\cap C^1(M)$ with $\int_M |\nabla_gu_0|^ndv_g=1$ and $\int_M u_0dv_g=0$ such that
\be\label{li-sci}\int_M e^{\alpha_n |u_0|^{n/(n-1)}}dv_g=\sup_{u\in W^{1,n}(M),\,\int_M|\nabla_gu|^ndv_g\leq 1,\,\int_Mudv_g=0}
\int_M e^{\alpha_n |u|^{n/(n-1)}}dv_g.\ee
Obviously (\ref{li-sci}) implies (\ref{Fontana}).
  In \cite{Yang-Tran}, the inequality (\ref{A-D}) was generalized to a closed Riemannian surface version, namely for any $\alpha$ with $0\leq\alpha<\lambda_1(\Sigma)=\inf_{\|u\|_2=1,\,\int_\Sigma udv_g=0}\|\nabla_gu\|_2^2$,
\be\label{Yang-Tr}\sup_{u\in W^{1,2}(\Sigma,g),\,\int_\Sigma|\nabla_g u|^2dv_g\leq 1,\,\int_\Sigma udv_g=0}\int_\Sigma e^{4\pi u^2(1+\alpha\|u\|_2^2)}dv_g<\infty;\ee
moreover, the supremum in (\ref{Yang-Tr}) can be attained for sufficiently small $\alpha$. However, in a recent work \cite{Yang-JDE-15},
an analog of (\ref{Tintarev}) was also established  on a
closed Riemannian surface, say for any $\alpha<\lambda_1(\Sigma)$,
\be\label{Yang-JDE-15}\sup_{u\in W^{1,2}(\Sigma,g),\,\int_\Sigma|\nabla_g u|^2dv_g-\alpha \int_\Sigma u^2dv_g\leq 1}
\int_\Sigma e^{4\pi u^2}dv_g<\infty.\ee
Moreover, the above supremum can be attained for any $\alpha<\lambda_1(\Sigma)$.
Further, this kind of inequalities involving higher order eigenvalues of the Laplace-Beltrami operator
has been studied.

  In this paper, our aim is to establish Trudinger-Moser inequalities for ``symmetric" functions and prove the existence of their extremal functions
  on a closed Riemannian surface with the action of a finite isometric group. They can be viewed as a ``combination" of (\ref{Mosr-even}) and
  (\ref{Yang-JDE-15}). We believe that such inequalities would play an important role in the study of prescribing Gaussian curvature problem and
  mean field equations. Before ending this introduction, we mention Mancini-Martinazzi \cite{Mancini-Martinazzi}, who studied the classical
  Trudinger-Moser inequality by estimating the energy of extremals for subcritical functionals.

\section{Notations and main results}

Let $(\Sigma,g)$ be a closed Riemannian surface and
$\textbf{G}=\{\sigma_1,\cdots,\sigma_N\}$ be
an isometric group acting on it, where $N$ is some positive integer. By definition, $\textbf{G}$ is a group and each $\sigma_i:\Sigma\ra\Sigma$
is an isometric map, particularly $\sigma_i^\ast g_x=g_{\sigma_i(x)}$ for all $x\in \Sigma$.
Let  $u:\Sigma\ra\mathbb{R}$ be a measurable function, we say that $u\in\mathscr{I}_{\textbf{G}}$ if $u$ is $\textbf{G}$-invariant,
namely $u(\sigma_i(x))=u(x)$ for any $1\leq i\leq N$ and almost every $x\in\Sigma$.
We denote $W^{1,2}(\Sigma,g)$
the closure of $C^\infty(\Sigma)$ under the norm
$$\|u\|_{W^{1,2}(\Sigma,g)}=\le(\int_\Sigma\le(|\nabla_gu|^2+u^2\ri)dv_g\ri)^{1/2},$$
where $\nabla_g$ and $dv_g$ stand for the gradient operator and the Riemannian volume element respectively.
Define a Hilbert space
\be\label{HG}\mathscr{H}_\textbf{G}= \le\{u\in W^{1,2}(\Sigma,g)\cap \mathscr{I}_\textbf{G}: \int_\Sigma udv_g=0\ri\}\ee
with an inner product
$$\langle u,v\rangle_{\mathscr{H}_{\textbf{G}}}=\int_\Sigma \langle\nabla_gu,\nabla_gv\rangle dv_g,$$
where $\langle\nabla_gu,\nabla_gv\rangle$ stands for the Riemannian inner product of $\nabla_gu$ and $\nabla_gv$.
Let $\Delta_g=-{\rm div}_g\nabla_g$ be the Laplace-Beltrami operator, and
\be\label{eigenvalue-1}\lambda_1^{\textbf{G}}=\inf_{u\in \mathscr{H}_{\textbf{G}},\,u\not\equiv 0}\f{\int_\Sigma|\nabla_gu|^2dv_g}
{\int_\Sigma u^2dv_g}\ee
be the first eigenvalue of $\Delta_g$ on the space $\mathscr{H}_{\textbf{G}}$. For any $x\in\Sigma$, we set $I(x)=\sharp \textbf{G}(x)$, where $\sharp A$
stands for the number of all distinct points in the set $A$, and $\textbf{G}(x)=\{\sigma_1(x),\cdots,\sigma_N(x)\}$.
 Let
\be\label{ell}\ell=\min_{x\in\Sigma}I(x).\ee
Clearly we have $1\leq\ell\leq N$ since $1\leq I(x)\leq N$ for all $x\in\Sigma$. As one will see, the best constant in the Trudinger-Moser inequality for ``symmetric"
functions would be $4\pi\ell$. Precisely we state the following theorem.
\begin{theorem}\label{Thm1} Let $(\Sigma,g)$ be a closed Riemannian surface and $\textbf{G}=\{\sigma_1,\cdots,\sigma_N\}$ be an isometric group
acting on it. Assume $\mathscr{H}_{\textbf{G}}$, $\lambda_1^{\textbf{G}}$ and $\ell$ are defined by (\ref{HG}), (\ref{eigenvalue-1}) and (\ref{ell}) respectively.
Then we have the following assertions:\\
$(i)$ For any $\alpha<\lambda_1^\textbf{G}$ and $\beta\leq 4\pi\ell$, there holds
\be\label{ineq-1}\sup_{u\in\mathscr{H}_\textbf{G},\,\int_\Sigma|\nabla_gu|^2dv_g-\alpha \int_\Sigma u^2dv_g\leq 1}\int_\Sigma e^{\beta u^2}dv_g<\infty;\ee
$(ii)$ If $\alpha<\lambda_1^\textbf{G}$ and $\beta>4\pi\ell$, or $\alpha\geq \lambda_1^\textbf{G}$ and $\beta>0$, then the supremum in (\ref{ineq-1}) is infinity;\\
$(iii)$ If $\alpha<\lambda_1^\textbf{G}$ and $\beta\leq 4\pi\ell$, then the supremum in (\ref{ineq-1}) can be attained by some function
$u_0\in\mathscr{H}_\textbf{G}\cap C^1(\Sigma,g)$ with $\int_\Sigma|\nabla_gu_0|^2dv_g-\alpha \int_\Sigma u_0^2dv_g=1$.
\end{theorem}

As in \cite{Yang-JDE-15}, we may consider the effect of higher order eigenvalues on the Trudinger-Moser inequality.
For this purpose, we define the first eigenfunction space with respect to $\lambda_1^\textbf{G}$ by
$$E_{\lambda_1^\textbf{G}}=\le\{u\in \mathscr{H}_\textbf{G}: \Delta_gu=\lambda_1^\textbf{G} u\ri\}.$$
 By an induction,  the $j$-th ($j\geq 2$) eigenvalue and
eigenfunction space will be defined as
\be\label{lamj}\lambda_{j}^\textbf{G}=\inf_{u\in\mathscr{H}_\textbf{G},\,u\in E_{j-1}^\perp,\,u\not\equiv 0}
\f{\int_\Sigma|\nabla_gu|^2dv_g}{\int_\Sigma u^2dv_g}\ee
and
$$E_{\lambda_j^\textbf{G}}=\le\{u\in E_{j-1}^\perp: \Delta_gu=\lambda_j^\textbf{G} u\ri\}$$
respectively,
where $E_{j-1}=E_{\lambda_1^\textbf{G}}\oplus\cdots\oplus E_{\lambda_{j-1}^\textbf{G}}$ and
\be\label{Ej}E_{j-1}^\perp=\le\{u\in \mathscr{H}_\textbf{G}: \int_\Sigma uvdv_g=0,\,\forall v\in E_{j-1}\ri\}.\ee
Then higher order eigenvalues of $\Delta_g$ affect the Trudinger-Moser inequality in the following way:
\begin{theorem}\label{Thm2} Let $(\Sigma,g)$ be a closed Riemannian surface and $\textbf{G}=\{\sigma_1,\cdots,\sigma_N\}$ be an isometric group
acting on it. Assume $\mathscr{H}_\textbf{G}$, $\ell$, $\lambda_j^\textbf{G}$ and $E_{j-1}^\perp$ are defined by (\ref{HG}), (\ref{ell}), (\ref{lamj}) and (\ref{Ej}) respectively, $j\geq 2$.\\
$(i)$ For any $\alpha<\lambda_j^\textbf{G}$ and $\beta\leq 4\pi\ell$, there holds
\be\label{ineq-2}\sup_{u\in E_{j-1}^\perp,\,\int_\Sigma|\nabla_gu|^2dv_g-\alpha \int_\Sigma u^2dv_g\leq 1}\int_\Sigma e^{\beta u^2}dv_g<\infty;\ee
$(ii)$ If $\alpha<\lambda_j^\textbf{G}$ and $\beta>4\pi\ell$, or $\alpha\geq \lambda_j^\textbf{G}$ and $\beta>0$, then the supremum in (\ref{ineq-2}) is infinity;\\
$(iii)$ For any $\alpha<\lambda_j^\textbf{G}$ and $\beta\leq 4\pi\ell$, the supremum in (\ref{ineq-2}) can be attained by some function
$u_0\in E_{j-1}^\perp\cap C^1(\Sigma,g)$ with $\int_\Sigma|\nabla_gu_0|^2dv_g-\alpha \int_\Sigma u_0^2dv_g=1$.
\end{theorem}

Let us give several examples for the finite isometric group $\textbf{G}$ acting on a closed Riemannian surface $(\Sigma,g)$.  $(a)$ If $\textbf{G}=\{Id\}$, where $Id$ denotes the identity map, then $\textbf{G}$ is a trivial isometric group action, and Theorems \ref{Thm1} and \ref{Thm2} are reduced to (\cite{Yang-JDE-15}, Theorems 3 and 4).
$(b)$ Let $(\mathbb{S}^2,g_0)$ be the standard 2-sphere given as in the introduction, and $\textsl{\textbf{G}}=\{Id, \sigma_0\}$, where $\sigma_0(x)=-x$ for any $x\in \mathbb{S}^2$. Then we have $\sharp \textbf{G}(x)=\sharp\{x,-x\}=2$ for any $x\in \mathbb{S}^2$, and thus $\ell=2$. Hence  Moser's inequality (\ref{Mosr-even})
for even functions  is a special case of our theorems. $(c)$ If $\textbf{G}$ has a fixed point, namely there exists some point
$p\in \Sigma$ such that $\sigma(p)=p$ for all $\sigma\in \textbf{G}$, then we have $\ell=\sharp\textbf{G}(p)=1$, and whence both of the best constants
in (\ref{ineq-1}) and (\ref{ineq-2}) are $4\pi$.

From now on, to simplify notations, we write
\be\label{norm-alpha}\|u\|_{1,\alpha}=\le(\int_\Sigma|\nabla_gu|^2dv_g-\alpha \int_\Sigma u^2dv_g\ri)^{1/2},\ee
provided that the right hand side of the above equality makes sense, say, if $\alpha<\lambda_1^\textbf{G}$ and
$u\in\mathscr{H}_\textbf{G}$, then $\|u\|_{1,\alpha}$ is well defined.
For the proof of Theorems \ref{Thm1} and \ref{Thm2}, we follow the lines of \cite{Yang-JDE-15} and thereby follow
closely \cite{Lijpde}. Pioneer works are due to Carleson-Chang \cite{C-C},
Ding-Jost-Li-Wang \cite{DJLW}, and Adimurthi-Struwe \cite{A-S}.
Since both of them are similar, we only give the outline of the proof of Theorem \ref{Thm1}. Firstly, we prove that the best constant in (\ref{ineq-1}) is $4\pi\ell$, which is based on Moser's original inequality and test function computations; Secondly,
 a direct method of variation shows that every subcritical Trudinger-Moser functional has a maximizer, namely for any $\epsilon>0$, there exists
 some $u_\epsilon\in\mathscr{H}_\textbf{G}$ with $\|u_\epsilon\|_{1,\alpha}=1$ satisfying
 $$\int_\Sigma e^{(4\pi\ell-\epsilon)u_\epsilon^2}dv_g=\sup_{u\in\mathscr{H}_\textbf{G},\,\|u\|_{1,\alpha}\leq 1}\int_\Sigma e^{(4\pi\ell-\epsilon)u^2}dv_g,$$
 where $\alpha<\lambda_1^\textbf{G}$ and $\|u\|_{1,\alpha}$ is defined as in (\ref{norm-alpha});
 Thirdly, we use blow-up analysis to show that if $\sup_{x\in \Sigma}|u_\epsilon|\ra\infty$ as $\epsilon\ra 0$, then
 $$\sup_{u\in\mathscr{H}_\textbf{G},\,\|u\|_{1,\alpha}\leq 1}\int_\Sigma e^{4\pi\ell u^2}dv_g\leq {\rm Vol}_g(\Sigma)+
\pi\ell e^{1+4\pi\ell A_{x_0}},$$
where $A_{x_0}$ is a constant related to certain Green function (see (\ref{ax00}) below); Finally, we construct a sequence of
functions $\phi_\epsilon\in\mathscr{H}_\textbf{G}$ with $\|\phi_\epsilon\|_{1,\alpha}\leq 1$ such that
$$\int_\Sigma e^{4\pi\ell \phi_\epsilon^2}dv_g>{\rm Vol}_g(\Sigma)+
\pi\ell e^{1+4\pi\ell A_{x_0}},$$
provided that $\epsilon>0$ is chosen sufficiently small. Combining the above two estimates, we get a contradiction,  which implies that
$u_\epsilon$ must be uniformly bounded. Then applying elliptic estimates to the equation of $u_\epsilon$, we get a desired extremal
function.

In the remaining part of this paper, we shall prove Theorems \ref{Thm1} and \ref{Thm2}.
Throughout this paper, we do not distinguish sequence and subsequence. Moreover we often denote various constants by the same $C$,
but the dependence of $C$ will be given only if necessary. Also we use  symbols $|O(R\epsilon)|\leq CR\epsilon$,  $o_\epsilon(1)\ra 0$ as $\epsilon\ra 0$,
$o_\delta(1)\ra 0$ as $\delta\ra 0$, and so on.

\section{Proof of Theorem \ref{Thm1}}\label{sec-2}

In this section, we shall prove Theorem \ref{Thm1}. In the first subsection, we show that the best
constant in (\ref{ineq-1}) is equal to $4\pi\ell$. The essential tools we use are subcritical Trudinger-Moser inequality
and Moser's sequence of functions. Also we prove $(ii)$ of Theorem \ref{Thm1}. In the second subsection, we consider the existence of
maximizers for subcritical Trudinger-Moser functionals and study their energy concentration phenomenon. In the third subsection,
assuming blow-up occurs, we derive an upper bound of the supremum in (\ref{ineq-1}), which obviously leads to
$(i)$ of Theorem \ref{Thm1}.  In the final subsection, we construct
a sequence of test functions to show that the upper bound we obtained in the third subsection is not really an upper bound. Therefore blow-up can not occur
and elliptic estimates lead to existence of extremal function. This concludes $(iii)$ of Theorem \ref{Thm1}.

\subsection{The best constant}\label{Subsec-1}
In view of (\ref{eigenvalue-1}), one can see that
$\lambda_1^\textbf{G}>0$ by using a direct method of variation. For any fixed $\alpha<\lambda_1^\textbf{G}$,
if $u\in\mathscr{H}_\textbf{G}$ satisfies $\|u\|_{1,\alpha}\leq 1$, then $\|\nabla_gu\|_2^2
\leq {\lambda_1^\textbf{G}}/{(\lambda_1^\textbf{G}-\alpha)}$. By Fontana's inequality (\ref{Fontana}),
there exists a positive constant $\beta_0$ depending only on
$\lambda_1^\textbf{G}$ and $\alpha$ such that
$$\sup_{u\in\mathscr{H}_\textbf{G},\,\|u\|_{1,\alpha}\leq 1}\int_\Sigma e^{\beta_0u^2}dv_g<\infty.$$
Now we define
\be\label{best-constant}\beta^\ast=\sup\le\{\beta:\sup_{u\in\mathscr{H}_\textbf{G},\,\|u\|_{1,\alpha}\leq 1}\int_\Sigma e^{\beta
u^2}dv_g<\infty
\ri\}.\ee

\begin{lemma}\label{lemma-1} Let $\ell$ and $\beta^\ast$ be defined as in (\ref{ell}) and (\ref{best-constant}) respectively. Then
$\beta^\ast=4\pi\ell$.
\end{lemma}
{\it Proof}. We divide the proof into two steps.

{\it Step 1. There holds $\beta^\ast\leq 4\pi\ell$.}

In view of (\ref{ell}), there exists some point $x_0\in\Sigma$ satisfying
$\ell=\sharp \textbf{G}(x_0)=\sharp\le\{\sigma_1(x_0),\cdots,\sigma_N(x_0)\ri\}$.
Without loss of generality, we assume that $\sigma_1=Id$ is the identity map, and that
$\textbf{G}(x_0)=\{\sigma_i(x_0)\}_{i=1}^\ell$.
Take
$$r_0=\f{1}{4}\min_{1\leq i<j\leq \ell}d_g(\sigma_i(x_0),\sigma_j(x_0)),$$
where $d_g(\sigma_i(x_0),\sigma_j(x_0))$ denotes the Riemannian distance between $\sigma_i(x_0)$ and $\sigma_j(x_0)$.
Since every $\sigma_i:\Sigma\ra\Sigma$ is an isometric map, we can see that for all $0<r\leq r_0$,
\be\label{Br}B_r(\sigma_i(x_0))=\sigma_i(B_r(x_0)), \,\,\, 1\leq i\leq\ell,\ee
where $B_r(x)$ stands for the geodesic ball centered at $x\in\Sigma$ with radius $r$.

Fixing $p\in\Sigma$, $k\in\mathbb{N}$ and $0<r\leq r_0$, we take a sequence of Moser functions by
\be\label{19}M_{p,k}=M_{p,k}(x,r)=\le\{
\begin{array}{lll}
\log k&{\rm when}&\rho\leq rk^{-1/4}\\[1.5ex]
4\log\f{r}{\rho}&{\rm when}&rk^{-1/4}<\rho\leq r\\[1.5ex]
0&{\rm when}& \rho>r,
\end{array}
\ri.\ee
where $\rho$ denotes the Riemannian distance between $x$ and $p$. Define
\be\label{test}\widetilde{M}_k=\widetilde{M}_k(x,r)=\le\{
\begin{array}{lll}
M_{\sigma_i(x_0),k}(x,r),& x\in B_{r_0}(\sigma_i(x_0)),\,\,1\leq i\leq\ell\\[1.5ex]
0, &x\in \Sigma\setminus\cup_{i=1}^\ell B_{r_0}(\sigma_i(x_0))
\end{array}
\ri.\ee
If $x\in B_{r_0}(\sigma_i(x_0))$ for some $i$, then it follows from (\ref{Br}) that
for any $j=1,\cdots,N$, $\sigma_j(x)\in B_r(\sigma_j(\sigma_i(x_0)))$ and
$d_g(\sigma_j(x),\sigma_j(\sigma_i(x_0)))=d_g(x,\sigma_i(x_0))$. In view of (\ref{19}) and (\ref{test}), one can easily check that
\be\label{invariant}\widetilde{M}_k(\sigma_j(x),r)=\widetilde{M}_k(x,r),\,\,\forall x\in B_{r_0}(\sigma_i(x_0)),\,\,1\leq i\leq \ell,\,\,1\leq j\leq N.\ee
If $x\in \Sigma\setminus\cup_{i=1}^\ell B_{r_0}(\sigma_i(x_0))$, then $\sigma_j(x)\in \Sigma\setminus\cup_{i=1}^\ell B_{r_0}(\sigma_i(x_0))$, and thus
$\widetilde{M}_k(\sigma_j(x),r)=0$ for $j=1,\cdots,N$. This together with (\ref{invariant}) leads to
\be\label{inv-t}\widetilde{M}_k(\sigma_j(x),r)=\widetilde{M}_k(x,r),\,\,\forall x\in\Sigma,\,\,1\leq j\leq N.\ee

A straightforward calculation shows
\bea\label{energy}
&& \int_\Sigma|\nabla_g\widetilde{M}_k|^2dv_g=(1+O(r))8\pi\ell \log k,\\
&& \int_\Sigma \widetilde{M}_k^mdv_g=O(1), \,\,m=1,2.\label{mean-v}
\eea
Denote
$\overline{\widetilde{M}}_k=\f{1}{{\rm Vol}_g(\Sigma)}\int_\Sigma \widetilde{M}_k dv_g$ and define
$$M_k^\ast=M_k^\ast(x,r)=\f{\widetilde{M}_k(x,r)-\overline{\widetilde{M}}_k}
{\|\widetilde{M}_k-\overline{\widetilde{M}}_k\|_{1,\alpha}}.$$
In view of (\ref{inv-t}), we have $M_k^\ast\in \mathscr{H}_\textbf{G}$. Note that $\|M_k^\ast\|_{1,\alpha}=1$.
 By (\ref{energy}) and (\ref{mean-v}),
$$\|\widetilde{M}_k-\overline{\widetilde{M}}_k\|_{1,\alpha}=(1+O(r))8\pi\ell \log k+O(1).$$
Hence we have for any $\beta_1>4\pi\ell$,
\bna
\int_{B_{r_0}(x_0)}e^{\beta_1 {M_k^\ast}^2}dv_g&\geq&\int_{B_{rk^{-1/4}}(x_0)}e^{\beta_1
\f{(\log k+O(1))^2}{(1+O(r))8\pi\ell\log k+O(1)}}dv_g\\
&=&e^{\f{\beta_1(1+o_k(1))\log k}{(1+O(r))8\pi\ell}}\pi r^2k^{-1/2}(1+o_k(1)).
\ena
Choosing $r>0$ sufficiently small and then passing to the limit $k\ra\infty$ in the above estimate, we
conclude
$$\int_{B_{r_0}(x_0)}e^{\beta_1 {M_k^\ast}^2}dv_g\ra \infty\quad{\rm as}\quad k\ra\infty.$$
Therefore $\beta^\ast\leq 4\pi\ell$.\\

{\it Step 2. There holds $\beta^\ast\geq 4\pi\ell$.}

Suppose $\beta^\ast<4\pi\ell$. Then for any $k\in\mathbb{N}$, there is a $u_k\in\mathscr{H}_\textbf{G}$ with $\|u_k\|_{1,\alpha}\leq 1$
such that
\be\label{over}\int_\Sigma e^{(\beta^\ast+k^{-1})u_k^2}dv_g\ra\infty\quad{\rm as}\quad k\ra\infty.\ee
Since $\alpha<\lambda_1^\textbf{G}$, we can see that $u_k$ is bounded in $W^{1,2}(\Sigma,g)$. Up to a subsequence, we can assume that
$u_k$ converges to some function $u_0$ weakly in $W^{1,2}(\Sigma,g)$, strongly in $L^q(\Sigma,g)$, $\forall q>1$, and
for almost every  $x\in\Sigma$. Clearly $u_0\in \mathscr{H}_\textbf{G}$ and $\|u_0\|_{1,\alpha}\leq 1$. We now {\it claim} that $u_0\equiv 0$.
For otherwise, we have
\be\label{less-1}\|u_k-u_0\|_{1,\alpha}^2\leq 1-\|u_0\|_{1,\alpha}^2+o_k(1)\leq 1-\f{1}{2}\|u_0\|_{1,\alpha}^2<1\ee
for sufficiently large $k$. Given any $\epsilon>0$. We calculate
\bea\nonumber
\int_\Sigma e^{(\beta^\ast+k^{-1})u_k^2}dv_g&\leq& \int_\Sigma e^{(\beta^\ast+k^{-1})(1+\epsilon)(u_k-u_0)^2+Cu_0^2}dv_g\\
&\leq& C\le(\int_\Sigma e^{(\beta^\ast+k^{-1})(1+2\epsilon)(u_k-u_0)^2}dv_g\ri)^{\f{1+\epsilon}{1+2\epsilon}},\label{contr}
\eea
where $C$ is a constant depending only on $u_0$, $\beta^\ast$ and $\epsilon$. In view of (\ref{less-1}), one can find
a small $\epsilon>0$ and a large integer $k_0$ such that when $k\geq k_0$, there holds
$$(\beta^\ast+k^{-1})(1+2\epsilon)\|u_k-u_0\|_{1,\alpha}^2\leq \beta^\ast \le(1-8^{-1}\|u_0\|_{1,\alpha}^2\ri).$$
This together with (\ref{contr}) leads to
$$\int_\Sigma e^{(\beta^\ast+k^{-1})u_k^2}dv_g\leq C,$$
contradicting (\ref{over}). This confirms our claim $u_0\equiv 0$.

For any fixed $x\in\Sigma$, we let $I=I(x)=\sharp \textbf{G}(x)$. Without loss of generality, we assume that
$\sigma_1=Id$ and that
$\textbf{G}(x)=\{\sigma_1(x),\cdots,\sigma_I(x)\}$. There exists sufficiently small $r_1>0$ such that
$\cap_{i=1}^IB_{r_1}(\sigma_i(x))=\varnothing$. Since $\sigma_i$'s are all isometric maps, if $0<r\leq r_1$, then we have
$$\int_{B_r(\sigma_i(x))}|\nabla_gu_k|^2dv_g=
\int_{B_r(x)}|\nabla_gu_k|^2dv_g,\quad\forall 1\leq i\leq I.$$
Noting that $I\geq \ell$, $\|u_k\|_{1,\alpha}\leq 1$ and $u_0\equiv 0$, we have for $0<r\leq r_1$,
\be\label{1/l}\int_{B_r(x)}|\nabla_gu_k|^2dv_g\leq \f{1}{\ell}+o_k(1).\ee
Let $\zeta\in C_0^1(B_r(x))$, $0\leq \zeta\leq 1$, $\zeta\equiv 1$ on $B_{r/2}(x)$ and $|\nabla_g\zeta|\leq \f{2}{r}$.
This together with (\ref{1/l}) and $u_0\equiv 0$ implies that $\zeta u_k\in W_0^{1,2}(B_r(x))$ and
\be\label{ell-1}\int_{B_r(x)}|\nabla_g(\zeta u_k)|^2dv_g\leq \f{1}{\ell}+o_k(1).\ee
Take a normal coordinate system $(B_r(x),\exp_x^{-1};\{y\})$, where $y=(y_1,y_2)\in\mathbb{B}_r(0)\subset\mathbb{R}^2$,
and $\exp_x:\mathbb{B}_{r}(0)\ra B_r(x)$ denotes the exponential map. Let $\psi_k(y)=(\zeta u_k)(\exp_x(y))$, $y\in \mathbb{B}_r(0)$. In view of (\ref{ell-1}),
one easily gets
\bea\nonumber\int_{\mathbb{B}_r(0)}|\nabla_{\mathbb{R}^2}\psi_k(y)|^2dy&=&(1+O(r))\int_{B_r(x)}|\nabla_g(\zeta u_k)|^2dv_g\\
&\leq& (1+O(r))\le(\f{1}{\ell}+o_k(1)\ri),\label{ell-2}\eea
where $\nabla_{\mathbb{R}^2}$ denotes the usual gradient operator in $\mathbb{R}^2$. Also there holds $\psi_k\in W_0^{1,2}(\mathbb{B}_r(0))$ since
$\zeta u_k\in W_0^{1,2}(B_r(x))$.
Hence, if $K\in\mathbb{N}$ is chosen sufficiently large and $r>0$ is chosen sufficiently small, it then follows from (\ref{ell-2}) and Moser's inequality (\ref{Moser-1}) that
\bea\nonumber\int_{B_{r/2}(x)}e^{(\beta^\ast+k^{-1})u_k^2}dv_g&\leq& \int_{B_r(x)}e^{(\beta^\ast+k^{-1})(\zeta u_k)^2}dv_g\\\nonumber
&=&(1+O(r))\int_{\mathbb{B}_r(0)}e^{(\beta^\ast+k^{-1})\psi_k^2}dy\\\label{30'}
&\leq& C\eea
for some constant $C$ and all $k\geq K$. Since $(\Sigma,g)$ is compact, there exists some constant $C$ such that for all
$k\geq K$,
$$\int_{\Sigma}e^{(\beta^\ast+k^{-1})u_k^2}dv_g\leq C.$$
This contradicts (\ref{over}) again. Hence $\beta^\ast\geq 4\pi\ell$.

We finish the proof of the lemma by combining Steps 1 and 2. $\hfill\Box$ \\

We now clarify  the proof of $(ii)$ of Theorem \ref{Thm1}, which is partially implied by Lemma \ref{lemma-1}.

{\it Proof of $(ii)$ of Theorem \ref{Thm1}.} If $\alpha<\lambda_1^\textbf{G}$ and $\beta>4\pi\ell$, then Step 1 of the proof of Lemma
\ref{lemma-1} gives the desired result. In the following, we assume $\alpha\geq \lambda_1^\textbf{G}$ and $\beta>0$. By a direct method of
variation, one can find a function $u_0\not\equiv 0$ satisfying $u_0\in \mathscr{H}_\textbf{G}\cap C^1(\Sigma)$
and
$$\int_\Sigma |\nabla_gu_0|^2dv_g=\lambda_1^\textbf{G}\int_\Sigma u_0^2dv_g.$$
For any $t\in \mathbb{R}$, we have $tu_0\in\mathscr{H}_\textbf{G}$ and
$$\int_\Sigma|\nabla_g (tu_0)|^2dv_g-\alpha\int_\Sigma (tu_0)^2dv_g\leq 0.$$
Moreover, there holds
$$\int_\Sigma e^{\beta(tu_0)^2}dv_g\ra\infty\quad{\rm as}\quad t\ra\infty.$$
Again this gives the desired result. $\hfill\Box$

\subsection{Maximizers for subcritical functionals}

Let $\alpha<\lambda_1^\textbf{G}$. As in (\cite{Yang-JDE-15}, page 3183),
by Lemma \ref{lemma-1} and a direct method of variation, we can prove that for any $0<\epsilon<4\pi\ell$, there exists some $u_\epsilon\in\mathscr{H}_\textbf{G}$
with $\|u_\epsilon\|_{1,\alpha}=1$ such that
\be\label{maximizer}\int_\Sigma e^{(4\pi\ell-\epsilon)u_\epsilon^2}dv_g=\sup_{u\in\mathscr{H}_\textbf{G},\,\|u\|_{1,\alpha}\leq 1}\int_\Sigma
e^{(4\pi\ell-\epsilon)u^2}dv_g.\ee
The Euler-Lagrange equation for the maximizer $u_\epsilon$ reads
\be\label{E-L}\le\{\begin{array}{lll}
\Delta_gu_\epsilon-\alpha u_\epsilon=\f{1}{\lambda_\epsilon}u_\epsilon e^{(4\pi\ell-\epsilon)u_\epsilon^2}
-\f{\mu_\epsilon}{\lambda_\epsilon}\\[1.5ex]
u_\epsilon\in\mathscr{H}_\textbf{G},\,\,\|u_\epsilon\|_{1,\alpha}=1\\[1.5ex]
\lambda_\epsilon=\int_\Sigma u_\epsilon^2e^{(4\pi\ell-\epsilon)u_\epsilon^2}dv_g\\[1.5ex]
\mu_\epsilon=\f{1}{{\rm Vol}_g(\Sigma)}\int_\Sigma u_\epsilon e^{(4\pi\ell-\epsilon)u_\epsilon^2}dv_g.
\end{array}\ri.\ee
Regularity theory implies that $u_\epsilon\in C^1(\Sigma,g)$. Using an argument of (\cite{Yang-JDE-15}, page 3184), one has
\be\label{posi}\liminf_{\epsilon\ra 0}\lambda_\epsilon>0,\quad |\mu_\epsilon|/\lambda_\epsilon\leq C.\ee
By (\ref{maximizer}), one can easily see that
\be\label{sup}\lim_{\epsilon\ra 0}\int_\Sigma e^{(4\pi\ell-\epsilon)u_\epsilon^2}dv_g=\sup_{u\in\mathscr{H}_G,\,\|u\|_{1,\alpha}\leq 1}
 \int_\Sigma e^{4\pi\ell u^2}dv_g.\ee
Note that we do {\it not} assume the supremum on the right hand side of (\ref{sup}) is finite.
If $|u_\epsilon|\leq C$, in view of (\ref{posi}), applying elliptic estimates to (\ref{E-L}), we obtain $u_\epsilon\ra u^\ast$ in $C^1(\Sigma,g)$,
which implies that $u^\ast\in\mathscr{H}_\textbf{G}$ and $\|u^\ast\|_{1,\alpha}=1$.
In view of (\ref{sup}), we know that $u^\ast$ is a desired extremal function. From now on, we assume
$c_\epsilon=\max_\Sigma|u_\epsilon|\ra +\infty$
as $\epsilon\ra 0$.
Noting that $-u_\epsilon$ also satisfies (\ref{maximizer}) and (\ref{E-L}), we may assume with no loss of generality that
\be\label{c-e}c_\epsilon=\max_\Sigma|u_\epsilon|=\max_\Sigma u_\epsilon=u_\epsilon(x_\epsilon)\ra+\infty\ee
 and that
 \be\label{x0}x_\epsilon\ra x_0\in\Sigma\quad{as}\quad \epsilon\ra 0.\ee
 To proceed, we need the following energy concentration phenomenon of $u_\epsilon$.
 \begin{lemma}\label{e-ll} Under the assumptions (\ref{c-e}) and (\ref{x0}), we have\\
 $(i)$ $u_\epsilon$ converges to $0$ weakly in $W^{1,2}(\Sigma,g)$, strongly in $L^2(\Sigma,g)$,
 and almost everywhere in $\Sigma$;\\
 $(ii)$  $I(x_0)=\sharp \textbf{G}(x_0)=\ell$; \\
 $(iii)$ $\lim\limits_{r\ra 0}\lim\limits_{\epsilon\ra 0}\int_{B_r(x_0)}|\nabla_gu_\epsilon|^2dv_g=1/\ell$.
 \end{lemma}
 {\it Proof.} $(i)$ Since $\alpha<\lambda_1^\textbf{G}$ and $\|u_\epsilon\|_{1,\alpha}=1$, $u_\epsilon$ is bounded in $W^{1,2}(\Sigma,g)$. Hence we may assume
 $u_\epsilon$ converges to $u_0$ weakly in $W^{1,2}(\Sigma,g)$, strongly in $L^2(\Sigma,g)$,
 and almost everywhere in $\Sigma$. If $u_0\not\equiv 0$, then
 $$\|u_\epsilon-u_0\|_{1,\alpha}^2=1-\|u_0\|_{1,\alpha}^2+o_\epsilon(1)\leq 1-\f{1}{2}\|u_0\|_{1,\alpha}^2,$$
 provided that $\epsilon$ is sufficiently small. It follows from Lemma \ref{lemma-1} that $e^{(4\pi\ell-\epsilon)u_\epsilon^2}$
 is bounded in $L^q(\Sigma,g)$ for some $q>1$. Then applying elliptic estimates to (\ref{E-L}), we have that
 $\|u_\epsilon\|_{L^\infty(\Sigma)}\leq C$, which contradicts (\ref{c-e}). Therefore $u_0\equiv 0$.

 $(ii)$ Since $\ell=\min_{x\in\Sigma} I(x)$, we have $I(x_0)\geq \ell$. Suppose $I=I(x_0)>\ell$.
 Using the same argument as we derived (\ref{1/l}), we have
 \be\label{ell-3}\int_{B_r(x_0)}|\nabla_gu_\epsilon|^2dv_g\leq \f{1}{I}+o_\epsilon(1),\ee
 provided that $r>0$ is chosen sufficiently small. Similar to (\ref{30'}), it follows from (\ref{ell-3}) and Moser's inequality (\ref{Moser-1}) that
  $$\int_{B_{r/2}(x_0)}e^{4\pi\ell p u_\epsilon^2}dv_g\leq C$$
  for some sufficiently small $r>0$ and some $p>1$, where $C$ is a constant depending only on $r$, $p$, $I$ and $\ell$.
  Applying elliptic estimates
 to (\ref{E-L}), we have that $u_\epsilon$ is uniformly bounded in $B_{r/4}(x_0)$. This contradicts (\ref{c-e}). Therefore $I(x_0)=\ell$.

 $(iii)$ By $(ii)$, there exists some $r_0>0$ such that $\|\nabla_gu_\epsilon\|_{L^2(B_{r_0}(x_0))}^2\leq \f{1}{\ell}+o_\epsilon(1)$. It
 follows that
 \be\label{lim-r}\lim_{r\ra 0}\lim_{\epsilon\ra 0}\int_{B_{r}(x_0)}|\nabla_gu_\epsilon|^2dv_g\leq \f{1}{\ell}.\ee
 We {\it claim} that the equality of (\ref{lim-r}) holds. For otherwise, there exist two positive constants $\nu$ and $r_1$ with $0<r_1<r_0$  such that
 $$\int_{B_{r_1}(x_0)}|\nabla_gu_\epsilon|^2dv_g< \f{1}{\ell}-\nu.$$
 Similarly as we did in the proof of $(ii)$, we have that $e^{(4\pi\ell-\epsilon)u_\epsilon^2}$
 is bounded in $L^q(B_{r_1/2}(x_0))$ for some $q>1$. Then applying elliptic estimates to (\ref{E-L}), we obtain that
 $u_\epsilon$ is uniformly bounded in $B_{r_1/4}(x_0)$, which contradicts (\ref{c-e}). This concludes our claim and $(iii)$ holds.
 $\hfill\Box$

 \subsection{Blow-up analysis}
 Set
 \be\label{scale}r_\epsilon=\f{\sqrt{\lambda_\epsilon}}{c_\epsilon}e^{-(2\pi\ell-\epsilon/2)c_\epsilon^2}.\ee
 For any $0<a<4\pi\ell$, by Lemma \ref{lemma-1}, the H\"older inequality and $(i)$ of Lemma \ref{e-ll}, one has
 $$\lambda_\epsilon=\int_\Sigma u_\epsilon^2 e^{(4\pi\ell-\epsilon)u_\epsilon^2}dv_g= e^{a c_\epsilon^2}
 \int_\Sigma u_\epsilon^2e^{(4\pi\ell-\epsilon-a)u_\epsilon^2}dv_g\leq e^{ac_\epsilon^2}o_\epsilon(1).$$
 It then follows that
 \be\label{r-0}r_\epsilon^2c_\epsilon^2e^{(4\pi\ell-\epsilon-a)c_\epsilon^2}=o_\epsilon(1).\ee
 In particular, $r_\epsilon\ra 0$ as $\epsilon\ra 0$.
 Let $0<\delta<\f{1}{2}i_g(\Sigma)$ be fixed, where $i_g(\Sigma)$ is the injectivity radius of $(\Sigma,g)$. For
 $y\in\mathbb{B}_{\delta r_\epsilon^{-1}}(0)\subset\mathbb{R}^2$,
 we define $\psi_\epsilon(y)=c_\epsilon^{-1}u_\epsilon(\exp_{x_\epsilon}(r_\epsilon y))$, $\varphi_\epsilon(y)=c_\epsilon(u_\epsilon(\exp_{x_\epsilon}(r_\epsilon y))-c_\epsilon)$ and
 $g_\epsilon(y)=(\exp^\ast_{x_\epsilon}g)(r_\epsilon y)$,
 where $\mathbb{B}_{\delta r_\epsilon^{-1}}(0)$ is the Euclidean ball of radius $\delta r_\epsilon^{-1}$ centered at $0$, and
 $\exp_{x_\epsilon}$ is the exponential map at $x_\epsilon$. Note that $g_\epsilon$ converges to $g_0$
 in $C^2_{\rm loc}(\mathbb{R}^2)$ as $\epsilon\ra 0$, where $g_0$ denotes the standard Euclidean metric. By (\ref{E-L}), we have on $\mathbb{B}_{\delta r_\epsilon^{-1}}(0)$,
 \bea
 \label{psi-eqn}
 &&\Delta_{g_\epsilon}\psi_\epsilon(y)=\alpha r_\epsilon^2\psi_\epsilon(y)+c_\epsilon^{-2}\psi_\epsilon(y)e^{(4\pi\ell-\epsilon)
 (u_\epsilon^2(\exp_{x_\epsilon}(r_\epsilon y))-c_\epsilon^2)}-r_\epsilon^2c_\epsilon^{-1}\f{\mu_\epsilon}{\lambda_\epsilon}
 \\
 \label{varphi-e}&&\Delta_{g_\epsilon}\varphi_\epsilon(y)=\alpha r_\epsilon^2c_\epsilon^2\psi_\epsilon(y)+\psi_\epsilon(y)
 e^{(4\pi\ell-\epsilon)(u_\epsilon^2(\exp_{x_\epsilon}(r_\epsilon y))-c_\epsilon^2)}-r_\epsilon^2c_\epsilon\f{\mu_\epsilon}{\lambda_\epsilon}.
 \eea
 In view of (\ref{r-0}), applying elliptic estimates to (\ref{psi-eqn}) and (\ref{varphi-e}) respectively, we have
 \be\label{psi-1}\psi_\epsilon \ra 1\quad{\rm in}\quad C^1_{\rm loc}(\mathbb{R}^2),\ee
 and
 \be\label{vp}\varphi_\epsilon\ra \varphi\quad {\rm in}\quad C^1_{\rm loc}(\mathbb{R}^2),\ee
 where $\varphi$ satisfies
 $$\le\{\begin{array}{lll}
 -\Delta_{\mathbb{R}^2}\varphi=e^{8\pi\ell \varphi}\quad{\rm in}\quad \mathbb{R}^2\\[1.5ex]
 \varphi(0)=0=\sup_{\mathbb{R}^2}\varphi\\[1.5ex]
 \int_{\mathbb{R}^2}e^{8\pi\ell \varphi(y)}dy<\infty.
 \end{array}\ri.$$
 By a result of Chen-Li \cite{C-L}, we have
 $$\varphi(y)=-\f{1}{4\pi\ell}\log(1+\pi\ell|y|^2),$$
 which leads to
 \be\label{whole}\int_{\mathbb{R}^2}e^{8\pi\ell\varphi(y)}dy=\f{1}{\ell}.\ee
 By (\ref{scale}), (\ref{psi-1}) and (\ref{vp}),  there holds for any $R>0$,
 \bna\int_{\mathbb{B}_R(0)}e^{4\pi\ell \varphi(y)}dy&=&\lim_{\epsilon\ra 0}
 \int_{\mathbb{B}_R(0)}e^{(4\pi\ell-\epsilon)(u_\epsilon^2(\exp_{x_\epsilon}(r_\epsilon y))-c_\epsilon^2)}dy\\
 &=&\lim_{\epsilon\ra 0}\f{c_\epsilon^2}{\lambda_\epsilon}\int_{B_{Rr_\epsilon}(x_\epsilon)}
 e^{(4\pi\ell-\epsilon)u_\epsilon^2}dv_g\\
 &=&\lim_{\epsilon\ra 0}\f{1}{\lambda_\epsilon}\int_{B_{Rr_\epsilon}(x_\epsilon)}u_\epsilon^2
 e^{(4\pi\ell-\epsilon)u_\epsilon^2}dv_g.\ena
 This together with (\ref{whole}) gives
 \be\label{er}\lim_{R\ra\infty}\lim_{\epsilon\ra 0}\f{1}{\lambda_\epsilon}\int_{B_{Rr_\epsilon}(x_\epsilon)}
 u_\epsilon^2
 e^{(4\pi\ell-\epsilon)u_\epsilon^2}dv_g=\f{1}{\ell}.\ee
 By $(ii)$ of Lemma \ref{e-ll} and (\ref{x0}), one has for all sufficiently small $\epsilon>0$,
 \be\label{dis-joint} \cap_{i=1}^\ell B_{Rr_\epsilon}(\sigma_i(x_\epsilon))=\varnothing.\ee
 Noting that $u_\epsilon\in \mathscr{H}_\textbf{G}$, we have
 $$\int_{B_{Rr_\epsilon}(\sigma_i(x_\epsilon))} u_\epsilon^2
 e^{(4\pi\ell-\epsilon)u_\epsilon^2}dv_g=\int_{B_{Rr_\epsilon}(x_\epsilon)} u_\epsilon^2
 e^{(4\pi\ell-\epsilon)u_\epsilon^2}dv_g,\quad  1\leq i\leq \ell.$$
 This together with (\ref{er}) and (\ref{dis-joint}) leads to
 \be\label{l-R}\lim_{R\ra\infty}\lim_{\epsilon\ra 0}\f{1}{\lambda_\epsilon}\int_{B_{Rr_\epsilon}(\sigma_i(x_\epsilon))}
 u_\epsilon^2
 e^{(4\pi\ell-\epsilon)u_\epsilon^2}dv_g=\f{1}{\ell},\quad 1\leq i\leq\ell.\ee
 By definition of $\lambda_\epsilon$ in (\ref{E-L}), we conclude from (\ref{l-R}) that
 \be\label{zero}\lim_{R\ra\infty}\lim_{\epsilon\ra 0}\f{1}{\lambda_\epsilon}
 \int_{\Sigma\setminus \cup_{i=1}^\ell B_{Rr_\epsilon}(\sigma_i(x_\epsilon))}
 u_\epsilon^2
 e^{(4\pi\ell-\epsilon)u_\epsilon^2}dv_g=0.\ee

 Similar to \cite{Lijpde,A-D}, $\forall 0<\beta<1$, we let $u_{\epsilon,\,\beta}=\min\{u_\epsilon,\beta c_\epsilon\}$.

 \begin{lemma}\label{ener-less}
 $\forall\, 0<\beta<1$, there holds
 $$\lim_{\epsilon\ra 0}\int_\Sigma |\nabla_gu_{\epsilon,\,\beta}|^2dv_g=\beta.$$
 \end{lemma}

 {\it Proof.} Multiplying (\ref{E-L}) by $u_{\epsilon,\,\beta}$, we have
 \bea\nonumber
 \int_\Sigma|\nabla_gu_{\epsilon,\,\beta}|^2dv_g&=&\int_\Sigma\nabla_gu_{\epsilon,\,\beta}\nabla u_\epsilon dv_g\\\nonumber
 &=&\f{1}{\lambda_\epsilon}\int_\Sigma u_{\epsilon,\,\beta}u_\epsilon e^{(4\pi\ell-\epsilon)u_\epsilon^2}dv_g+
 \alpha\int_\Sigma u_{\epsilon,\,\beta}u_\epsilon dv_g-\f{\mu_\epsilon}{\lambda_\epsilon}\int_\Sigma u_{\epsilon,\,\beta}dv_g\\\nonumber
 &=&\f{1}{\lambda_\epsilon}\sum_{i=1}^\ell\int_{B_{Rr_\epsilon}(\sigma_i(x_\epsilon))}u_{\epsilon,\,\beta}u_\epsilon
 e^{(4\pi\ell-\epsilon)u_\epsilon^2}dv_g\\\label{na-beta}
 &&+\f{1}{\lambda_\epsilon}\int_{\Sigma\setminus \cup_{i=1}^\ell
 B_{Rr_\epsilon}(\sigma_i(x_\epsilon))}u_{\epsilon,\,\beta}u_\epsilon
 e^{(4\pi\ell-\epsilon)u_\epsilon^2}dv_g+o_\epsilon(1).
 \eea
 Note that $0\leq u_{\epsilon,\,\beta}u_\epsilon\leq u_\epsilon^2$ on $\Sigma$, and
 $u_{\epsilon,\,\beta}=\beta (1+o_\epsilon(1))u_\epsilon$ on $B_{Rr_\epsilon}(\sigma_i(x_\epsilon))$ for $1\leq i\leq\ell$.
 In view of (\ref{er}), (\ref{zero}) and (\ref{na-beta}), letting $\epsilon\ra 0$ first
 and then $R\ra\infty$, we conclude the lemma. $\hfill\Box$

 \begin{lemma}\label{l-c}
 There holds $\liminf_{\epsilon\ra 0} \lambda_\epsilon/c_\epsilon^2>0$.
 \end{lemma}

 {\it Proof.} Let $0<\beta<1$. In view of Lemma \ref{ener-less}, we have by using the H\"older inequality
 $$\int_{u_\epsilon\leq \beta c_\epsilon}u_\epsilon^2e^{(4\pi\ell-\epsilon)u_\epsilon^2}dv_g\leq \int_\Sigma
 u_\epsilon^2 e^{(4\pi\ell-\epsilon)u_{\epsilon,\,\beta}^2}dv_g=o_\epsilon(1).$$
 Similarly
 \bea\nonumber
 \f{\lambda_\epsilon}{c_\epsilon^2}&\geq& \beta^2\int_{u_\epsilon>\beta c_\epsilon}e^{(4\pi\ell-\epsilon)u_{\epsilon}^2}dv_g
 +o_\epsilon(1)\\\nonumber
 &\geq&\beta^2\le(\int_\Sigma e^{(4\pi\ell-\epsilon)u_{\epsilon}^2}dv_g-\int_\Sigma e^{(4\pi\ell-\epsilon)u_{\epsilon,\,\beta}^2}dv_g\ri)
 +o_\epsilon(1)\\\label{geq}
 &=&\beta^2\int_\Sigma (e^{(4\pi\ell-\epsilon)u_{\epsilon}^2}-1)dv_g+o_\epsilon(1).
 \eea
 This together with (\ref{sup}) ends the proof of the lemma. $\hfill\Box$

 \begin{lemma}\label{converge-Green} For any $1<q<2$, we have
 $c_\epsilon u_\epsilon$ converges to $G$ weakly in $W^{1,q}(\Sigma,g)$, strongly in $L^{2q/(2-q)}(\Sigma)$,
 and almost everywhere in $\Sigma$, where $G$ is a Green function satisfying
 \be\label{Green}\le\{\begin{array}{lll}
 \Delta_gG-\alpha G=\f{1}{\ell}\sum_{i=1}^\ell \delta_{\sigma_i(x_0)}-\f{1}{{\rm Vol}_g(\Sigma)}\\[1.5ex]
 \int_\Sigma Gdv_g=0\\[1.5ex]
 G(\sigma_i(x))=G(x),\, x\in \Sigma\setminus \{\sigma_j(x_0)\}_{j=1}^\ell,\, 1\leq i\leq\ell.
 \end{array}\ri.\ee
 \end{lemma}

 {\it Proof.} By (\ref{E-L}),
 \be\label{cu}\Delta_g(c_\epsilon u_\epsilon)-\alpha (c_\epsilon u_\epsilon)=h_\epsilon=\f{1}{\lambda_\epsilon}c_\epsilon u_\epsilon
 e^{(4\pi\ell-\epsilon)u_\epsilon^2}-\f{c_\epsilon \mu_\epsilon}{\lambda_\epsilon}.\ee
 It follows from Lemmas \ref{ener-less} and \ref{l-c} that for any $0<\beta<1$,
 \bna
 \int_\Sigma \f{c_\epsilon}{\lambda_\epsilon}|u_\epsilon|e^{(4\pi\ell-\epsilon)u_\epsilon^2}dv_g&=&
 \f{c_\epsilon}{\lambda_\epsilon}\int_{u_\epsilon\leq \beta c_\epsilon}|u_\epsilon| e^{(4\pi\ell-\epsilon)u_\epsilon^2}dv_g+
 \f{c_\epsilon}{\lambda_\epsilon}\int_{u_\epsilon> \beta c_\epsilon} u_\epsilon e^{(4\pi\ell-\epsilon)u_\epsilon^2}dv_g\\
 &\leq& \f{c_\epsilon}{\lambda_\epsilon}\int_\Sigma |u_\epsilon| e^{(4\pi\ell-\epsilon)u_{\epsilon,\,\beta}^2}dv_g+\f{1}{\beta}\\
 &\leq&\f{1}{\beta}+o_\epsilon(1),
 \ena
 and that
 \bna
 \f{c_\epsilon|\mu_\epsilon|}{\lambda_\epsilon}&\leq &\f{1}{{\rm Vol}_g(\Sigma)}\f{c_\epsilon}{\lambda_\epsilon}\int_{u_\epsilon\leq \beta c_\epsilon}
 |u_\epsilon| e^{(4\pi\ell-\epsilon)u_\epsilon^2}dv_g+\f{1}{{\rm Vol}_g(\Sigma)}\f{c_\epsilon}{\lambda_\epsilon}\int_{u_\epsilon> \beta c_\epsilon}
 u_\epsilon e^{(4\pi\ell-\epsilon)u_\epsilon^2}dv_g\\&\leq& \f{1}{{\rm Vol}_g(\Sigma)}\f{1}{\beta}+o_\epsilon(1).
 \ena
 Hence $h_\epsilon$ is bounded in $L^1(\Sigma,g)$. Then by (\cite{Yang-Zhu-SCM}, Lemma 2.11), we have $c_\epsilon u_\epsilon$
 is bounded in $W^{1,q}(\Sigma,g)$ for any $1<q<2$. Up to a subsequence, for any $1<q<2$ and $1<s\leq 2q/(2-q)$, $c_\epsilon u_\epsilon$
 converges to $G$ weakly in $W^{1,q}(\Sigma)$, strongly in $L^s(\Sigma,g)$, and almost everywhere in $\Sigma$.

 We calculate
 \be\label{est-1}\int_{u_\epsilon\leq \beta c_\epsilon}\f{c_\epsilon}{\lambda_\epsilon}u_\epsilon e^{(4\pi\ell-\epsilon)u_\epsilon^2}dv_g=o_\epsilon(1),\ee
 \be\label{est-2}
 \int_{\{u_\epsilon> \beta c_\epsilon\}\setminus \cup_{i=1}^\ell B_{Rr_\epsilon}(\sigma_i(x_\epsilon))}\f{c_\epsilon}{\lambda_\epsilon}u_\epsilon e^{(4\pi\ell-\epsilon)u_\epsilon^2}dv_g
 \leq \f{1}{\beta}\f{1}{\lambda_\epsilon}\int_{\Sigma\setminus \cup_{i=1}^\ell B_{Rr_\epsilon}(\sigma_i(x_\epsilon))}
 u_\epsilon^2 e^{(4\pi\ell-\epsilon)u_\epsilon^2}dv_g=o(1),
 \ee
 \be\label{est-3}\int_{B_{Rr_\epsilon}(\sigma_i(x_\epsilon))}\f{c_\epsilon}{\lambda_\epsilon}u_\epsilon e^{(4\pi\ell-\epsilon)u_\epsilon^2}dv_g
 =\f{1+o_\epsilon(1)}{\lambda_\epsilon}\int_{B_{Rr_\epsilon}(\sigma_i(x_\epsilon))}u_\epsilon^2 e^{(4\pi\ell-\epsilon)u_\epsilon^2}dv_g
 =\f{1}{\ell}+o(1),\, 1\leq i\leq\ell,\ee
 where $o(1)\ra 0$ as $\epsilon\ra 0$ first and then $R\ra \infty$.
 Integrating the equation (\ref{cu}), we have by combining (\ref{est-1})-(\ref{est-3}),
 $$\f{c_\epsilon\mu_\epsilon}{\lambda_\epsilon}{\rm Vol}_g(\Sigma)=\int_\Sigma\f{c_\epsilon}{\lambda_\epsilon}u_\epsilon e^{(4\pi\ell-\epsilon)u_\epsilon^2}dv_g
 =1+o_\epsilon(1).$$
 In view of (\ref{est-1})-(\ref{est-3}) again, testing the equation (\ref{cu}) by $\phi\in C^2(\Sigma)$ and passing to the limit $\epsilon\ra 0$, we have
 $$\int_\Sigma G\Delta_g\phi dv_g-\alpha\int_\Sigma G\phi dv_g=\f{1}{\ell}\sum_{i=1}^\ell \phi(\sigma_i(x_0))-\f{1}{{\rm Vol}_g(\Sigma)}\int_\Sigma
 \phi dv_g.$$
 Since $c_\epsilon u_\epsilon\in \mathscr{H}_\textbf{G}$, we have $\int_\Sigma Gdv_g=0$ and
 $G(\sigma_i(x))=G(x)$ for all $x\in\Sigma\setminus\{\sigma_1(x_0),\cdots,\sigma_\ell(x_0)\}$ and all $1\leq i\leq \ell$. $\hfill\Box$\\

 Let $${\psi}(x)=G(x)+\f{1}{2\pi\ell}\sum_{i=1}^\ell\log d_g(\sigma_i(x_0),x).$$
 It follows from (\ref{Green}) that the distributional Laplacian of $\psi$ belongs to $L^s(\Sigma,g)$ for some $s>2$. Then we have by
  elliptic estimates that $\psi\in C^1(\Sigma,g)$. Let $r_0=\f{1}{4}\min_{1\leq i<j\leq \ell}d_g(\sigma_i(x_0),\sigma_j(x_0))$.
  For $x\in B_{r_0}(x_0)$, the Green function $G$ can be decomposed as
  \be\label{G-representation}G(x)=-\f{1}{2\pi\ell}\log d_g(x,x_0)+A_{x_0}+\widetilde{\psi}(x),\ee
  where $\widetilde{\psi}\in C^1(\overline{B_{r_0}(x_0)})$, $\widetilde{\psi}(x_0)=0$ and
  \be\label{ax00}A_{x_0}=\lim_{x\ra x_0}\le(G(x)+\f{1}{2\pi\ell}\log d_g(x,x_0)\ri)=\lim_{x\ra x_0}\le(\psi(x)-
  \f{1}{2\pi\ell}\sum_{i=2}^\ell\log d_g(\sigma_i(x_0),x)\ri).\ee
   By (\ref{Green}), we have
 \bna
 \int_{\Sigma\setminus\cup_{i=1}^\ell B_\delta(\sigma_i(x_0))}|\nabla_gG|^2dv_g&=&\alpha
 \int_{\Sigma\setminus\cup_{i=1}^\ell B_\delta(\sigma_i(x_0))}G^2dv_g-\int_{\cup_{i=1}^\ell \p B_\delta(\sigma_i(x_0))}
 G\f{\p G}{\p\nu}d\sigma\\&&\quad-\f{1}{{\rm Vol}_g(\Sigma)}
 \int_{\Sigma\setminus\cup_{i=1}^\ell B_\delta(\sigma_i(x_0))}Gdv_g\\
 &=&-\f{1}{2\pi\ell}\log\delta+A_{x_0}+\alpha\int_\Sigma G^2dv_g+o_\delta(1).
 \ena
 Hence
 $$\int_{\Sigma\setminus\cup_{i=1}^\ell B_\delta(\sigma_i(x_0))}|\nabla_gu_\epsilon|^2dv_g
 =\f{1}{c_\epsilon^2}\le(-\f{1}{2\pi\ell}\log\delta+A_{x_0}+\alpha\int_\Sigma G^2dv_g+o_\delta(1)+o_\epsilon(1)\ri).$$
 It follows that
 \bna
 \int_{\cup_{i=1}^\ell B_\delta(\sigma_i(x_0))}|\nabla_gu_\epsilon|^2dv_g&=&1+\alpha\int_\Sigma u_\epsilon^2dv_g-
 \int_{\Sigma\setminus\cup_{i=1}^\ell B_\delta(\sigma_i(x_0))}|\nabla_gu_\epsilon|^2dv_g\\
 &=&1-\f{1}{c_\epsilon^2}\le(-\f{1}{2\pi\ell}\log\delta+A_{x_0}+o_\delta(1)+o_\epsilon(1)\ri).
 \ena
 Let $s_\epsilon=\sup_{\p B_\delta(x_0)}u_\epsilon$ and $\widetilde{u}_\epsilon=(u_\epsilon-s_\epsilon)^+$. Then $\widetilde{u}_\epsilon
 \in W_0^{1,2}(B_\delta(x_0))$, and satisfies
 $$\int_{\cup_{i=1}^\ell B_\delta(\sigma_i(x_0))}|\nabla_g\widetilde{u}_\epsilon|^2dv_g
 \leq\tau_\epsilon=1-\f{1}{c_\epsilon^2}\le(-\f{1}{2\pi\ell}\log\delta+A_{x_0}+o_\delta(1)+o_\epsilon(1)\ri)$$

 Now we choose an isothermal coordinate system $(U,\phi;\{x^1,x^2\})$ near $x_0$ such that
          $B_{2\delta}(x_0)\subset U$,
          $\phi(x_0)=0$, and the metric $g=e^h(d{x^1}^2+d{x^2}^2)$ for some function $h\in C^1(\phi(U))$
          with $h(0)=0$. Clearly, for any $\delta>0$, there exists some $c(\delta)>0$ with $c(\delta)\ra 0$ as
          $\delta\ra 0$ such that $\sqrt{g}\leq 1+c(\delta)$ and $\phi(B_\delta(p))\subset
          \mathbb{B}_{\delta(1+c(\delta))}(0)\subset\mathbb{R}^2$.
          Noting that $\widetilde{u}_\epsilon=0$ outside $B_\delta(p)$ for sufficiently small $\delta$, we have
          $$\int_{\mathbb{B}_{\delta(1+c(\delta))}(0)}
          |\nabla_{\mathbb{R}^2}(\widetilde{u}_\epsilon\circ \phi^{-1})|^2dx=\int_{\phi^{-1}(\mathbb{B}_{\delta(1+c(\delta))}(0))}
          |\nabla_g\widetilde{u}_\epsilon|^2dv_g=\int_{B_\delta(x_0)}|\nabla_g\widetilde{u}_\epsilon|dv_g\leq \f{\tau_\epsilon}{\ell}.
          $$
          This together with a result of Carleson-Chang \cite{C-C}  leads to
          \bea\nonumber
          \limsup_{\epsilon\ra 0}\int_{B_{\delta}(p)}(e^{4\pi\ell
          \widetilde{u}_\epsilon^2/\tau_\epsilon}-1)dv_g&\leq&\limsup_{\epsilon\ra 0}\,(1+c(\delta))\int_{\mathbb{B}_{\delta(1+c(\delta))}(0)}(e^{4\pi\ell
          (\widetilde{u}_\epsilon\circ\phi^{-1})^2/\tau_\epsilon}-1)dx\\\label{B-delta}
          &\leq&
          \pi\delta^2(1+c(\delta))^3e.
          \eea

          Note that $|u_\epsilon|\leq c_\epsilon$ and $u_\epsilon/c_\epsilon=1+o_\epsilon(1)$ on the geodesic ball
          $B_{Rr_\epsilon}(x_\epsilon)\subset{\Sigma}$. We estimate on $B_{Rr_\epsilon}(x_\epsilon)$,
          \bna
          (4\pi\ell-\epsilon)u_\epsilon^2&\leq&4\pi\ell(\widetilde{u}_\epsilon+s_\epsilon)^2\\
          &\leq&
          4\pi\ell\widetilde{u}_\epsilon^2+8\pi
          \ell s_\epsilon\widetilde{u}_\epsilon+o_\epsilon(1)\\
          &\leq&4\pi\ell\widetilde{u}_\epsilon^2-4\log\delta+8\pi\ell
          A_{x_0}+o(1)\\
          &\leq&4\pi\ell\widetilde{u}_\epsilon^2/\tau_\epsilon-2\log\delta+4\pi\ell
          A_{x_0}+o(1).
          \ena
          Therefore
          \bea\nonumber
          \int_{B_{Rr_\epsilon}(x_\epsilon)}e^{(4\pi\ell-\epsilon)
          u_\epsilon^2}dv_g&\leq& \delta^{-2}e^{4\pi\ell
          A_{x_0}+o(1)}\int_{B_{Rr_\epsilon}(x_\epsilon)}e^{4\pi\ell\widetilde{u}_\epsilon^2/\tau_\epsilon}dv_g\\
          \nonumber
          &=&\delta^{-2}e^{4\pi\ell
          A_{x_0}+o(1)}\int_{B_{Rr_\epsilon}(x_\epsilon)}(e^{4\pi\ell\widetilde{u}_\epsilon^2/\tau_\epsilon}-1)dv_g+o(1)\\
          \label{ii}
          &\leq&\delta^{-2}e^{4\pi\ell
          A_{x_0}+o(1)}\int_{B_\delta(x_0)}(e^{4\pi\ell\widetilde{u}_\epsilon^2/\tau_\epsilon}-1)dv_g+o(1),
          \eea
          where $o(1)\ra 0$ as $\epsilon\ra 0$ first and then $\delta \ra 0$.
          Combining (\ref{B-delta}) with (\ref{ii}), letting $\epsilon\ra 0$ first, and then letting $\delta\ra 0$,
          we conclude
          $$
          \limsup_{\epsilon\ra 0}\int_{B_{Rr_\epsilon}(x_\epsilon)}e^{(4\pi\ell-\epsilon)
          u_\epsilon^2}dv_g\leq \pi e^{1+4\pi\ell A_{x_0}}.
          $$
          Therefore
          \be\label{upper}\limsup_{\epsilon\ra 0}\int_{\cup_{i=1}^\ell B_{Rr_\epsilon}(\sigma_i(x_\epsilon))}e^{(4\pi\ell-\epsilon)
          u_\epsilon^2}dv_g\leq \pi\ell e^{1+4\pi\ell A_{x_0}}.\ee

\begin{proposition}\label{up-bd}
Under the assumptions $(\ref{c-e})$ and $(\ref{x0})$, there holds
$$\sup_{u\in\mathscr{H}_\textbf{G},\,\|u\|_{1,\alpha}\leq 1}\int_\Sigma e^{4\pi\ell u^2}dv_g=\lim_{\epsilon\ra 0}\int_\Sigma e^{(4\pi\ell-\epsilon)u_\epsilon^2}dv_g
\leq {\rm Vol}_g(\Sigma)+
\pi\ell e^{1+4\pi\ell A_{x_0}}.$$
\end{proposition}

{\it Proof.} We calculate
\bna
\int_{B_{Rr_\epsilon}(x_\epsilon)}e^{(4\pi\ell-\epsilon)u_\epsilon^2}dv_g&=&(1+o_\epsilon(1))
\int_{\mathbb{B}_R(0)}e^{(4\pi\ell-\epsilon)u_\epsilon^2(\exp_{x_\epsilon}(r_\epsilon y))}r_\epsilon^2dy\\
&=&(1+o_\epsilon(1))\f{\lambda_\epsilon}{c_\epsilon^2}\le(\int_{\mathbb{B}_R(0)}e^{8\pi\ell\varphi(y)}dy+o_\epsilon(1)\ri).
\ena
In view of (\ref{whole}) and (\ref{upper}),
$$\lim_{R\ra\infty}\lim_{\epsilon\ra 0}\int_{B_{Rr_\epsilon}(x_\epsilon)}e^{(4\pi\ell-\epsilon)u_\epsilon^2}dv_g=
\f{1}{\ell}\lim_{\epsilon\ra 0}\f{\lambda_\epsilon}{c_\epsilon^2}.$$
Hence
\be\label{c-lambda}\lim_{R\ra\infty}\lim_{\epsilon\ra 0}\int_{\cup_{i=1}^\ell B_{Rr_\epsilon}(\sigma_i(x_\epsilon))}
e^{(4\pi\ell-\epsilon)u_\epsilon^2}dv_g= \lim_{\epsilon\ra 0}\f{\lambda_\epsilon}{c_\epsilon^2}.\ee
By (\ref{geq}), we have
$$\lim_{\epsilon\ra 0}\int_\Sigma (e^{(4\pi\ell-\epsilon)u_\epsilon^2}-1)dv_g\leq \f{1}{\beta^2}\lim_{\epsilon\ra 0}
\f{\lambda_\epsilon}{c_\epsilon^2},\quad\forall 0<\beta<1.$$
Letting $\beta\ra 1$, we obtain
$$\lim_{\epsilon\ra 0}\int_\Sigma (e^{(4\pi\ell-\epsilon)u_\epsilon^2}-1)dv_g\leq \lim_{\epsilon\ra 0}
\f{\lambda_\epsilon}{c_\epsilon^2}.$$
This together with (\ref{upper}) and (\ref{c-lambda}) completes the proof of the proposition. $\hfill\Box$

\subsection{Test function computation}\label{Sec2.4}

In this subsection, we shall complete the proof of $(iii)$ of Theorem \ref{Thm1}. Let $\alpha<\lambda_1^{\textbf{G}}$ be fixed and
$\ell$ be an integer defined as in (\ref{ell}). In particular,
 we shall construct a function sequence
  $\phi_\epsilon$ satisfying $\phi_\epsilon\in\mathscr{H}_\textbf{G}$,
   \be\label{nrm}\int_\Sigma|\nabla_g\phi_\epsilon|^2dv_g-\alpha\int_\Sigma \phi_\epsilon^2dv_g=1\ee
   and
   \be\label{gr}\int_\Sigma e^{4\pi\ell \phi_\epsilon^2}dv_g>{\rm vol}_g(\Sigma)+\pi\ell e^{1+4\pi\ell A_{x_0}}\ee
   for sufficiently small $\epsilon>0$, where $x_0$ and $A_{x_0}$ are defined as in (\ref{x0}) and (\ref{ax00}) respectively.
   If there exists such a sequence $\phi_\epsilon$, then we have by Proposition \ref{up-bd}
    that $c_\epsilon$ must be bounded. Applying elliptic estimates
   to (\ref{E-L}), we conclude the existence of the desired extremal function.\\

   To do this, we define a sequence of functions by
   \be\label{b-e}\textsf{b}_\epsilon(x)=\le\{\begin{array}{lll}c+\f{-\f{1}{4\pi\ell}\log(1+\pi\ell\f{r^2}{\epsilon^2})+B}{c},&
   x\in B_{R\epsilon}(x_0)\\[1.5ex]\f{G-\zeta \widetilde{\psi}}{c}, & x\in B_{2R\epsilon}(x_0)\setminus
     B_{R\epsilon}(x_0),\end{array}\ri.\ee
   where $\widetilde{\psi}$ is defined as in (\ref{G-representation}), $\zeta\in C_0^\infty(B_{2R\epsilon}(x_0))$ satisfies that $\zeta\equiv 1$ on $B_{R\epsilon}(x_0)$ and
     $\|\nabla_g \zeta\|_{L^\infty}
     =O(1/(R\epsilon))$, $r=r(x)={\rm dist}_g(x,x_0)$, $R=-\log\epsilon$, $B$ and $c$ are constants depending only on $\epsilon$ to be determined later.
   Define another sequence of functions
   \be\label{eta-e}\eta_\epsilon(x)=\le\{\begin{array}{lll}
   \textsf{b}_\epsilon(x),&x\in B_{2R\epsilon}(x_0)\\[1.5ex]
   \textsf{b}_\epsilon(\sigma_i^{-1}(x)),&x\in B_{2R\epsilon}(\sigma_i(x_0)),\,2\leq i\leq\ell\\[1.5ex]
   \f{G}{c},&x\in \Sigma\setminus\cup_{i=1}^\ell B_{2R\epsilon}(\sigma_i(x_0)).
   \end{array}\ri.\ee
   Noting that $G(\sigma_i(x))=G(x)$ for all $x\in \Sigma\setminus\{\sigma_1(x_0),\cdots,\sigma_\ell(x_0)\}$,
   one can easily check that
   \be\label{h-g}\eta_\epsilon(\sigma_i(x))=\eta_\epsilon(x),\,\,\, \forall x\in\Sigma,\,\, \forall 1\leq i\leq\ell.\ee
     In view of (\ref{b-e}) and (\ref{eta-e}), in order to ensure that $\eta_\epsilon\in W^{1,2}(\Sigma,g),$ we set
     $$
     c+\f{1}{c}\le(-\f{1}{4\pi\ell}\log(1+\pi\ell R^2)+B\ri)
     =\f{1}{c}\le(-\f{1}{2\pi\ell}\log (R\epsilon)+A_{x_0}\ri),
     $$
     which gives
     \be\label{2pic2}
     2\pi \ell c^2=-\log\epsilon-2\pi\ell B+2\pi\ell A_{x_0}+\f{1}{2}\log (\pi\ell)
     +O(\f{1}{R^2}).
     \ee
     Noting that $\int_\Sigma Gdv_g=0$, we have
     \bea\nonumber
     \int_{\Sigma\setminus \cup_{i=1}^\ell B_{R\epsilon}(\sigma_i(x_0))}|\nabla_g G|^2dv_g&=&\int_{\Sigma\setminus \cup_{i=1}^\ell B_{R\epsilon}(\sigma_i(x_0))}G\Delta_gGdv_g-\int_{\cup_{i=1}^\ell\p B_{R\epsilon}(\sigma_i(x_0))}G\f{\p G}{\p\nu}d\sigma\\ \nonumber
     &=&\alpha\int_{\Sigma\setminus \cup_{i=1}^\ell B_{R\epsilon}(\sigma_i(x_0))}G^2dv_g-\f{1}{{\rm Vol}_g(\Sigma)}\int_{\Sigma\setminus \cup_{i=1}^\ell B_{R\epsilon}(\sigma_i(x_0))}Gdv_g\\
     \nonumber &&-\sum_{i=1}^\ell\int_{\p B_{R\epsilon}(\sigma_i(x_0))}G\f{\p G}{\p\nu}d\sigma\\\label{nabl-G}
     &=&-\f{1}{2\pi\ell}\log(R\epsilon)+\alpha\int_\Sigma G^2dv_g+A_{x_0}+O(R\epsilon\log(R\epsilon)).
     \eea
     Since $\widetilde{\psi}\in C^1(\Sigma,g)$ and $\widetilde{\psi}(x_0)=0$, we have
     \bea\label{nabl-eta}&&\int_{B_{2R\epsilon}(x_0)\setminus B_{R\epsilon}(x_0)}|\nabla_g\zeta|^2\widetilde{\psi}^2dv_g=O((R\epsilon)^2),\\
     \label{Geta}&&\int_{B_{2R\epsilon}(x_0)\setminus B_{R\epsilon}(x_0)}\nabla_g G\nabla_g\zeta\widetilde{\psi} dv_g=O(R\epsilon),\\
     \label{rep}&&\int_{B_{R\epsilon}(x_0)}|\nabla_g\eta_\epsilon|^2dv_g=\f{1}{\ell^2c^2}\le(\f{1}{2\pi}\log R+\f{\log(\pi\ell)}
     {4\pi}-\f{1}{4\pi}+O(\f{1}{R^2})\ri).\eea
     Combining (\ref{nabl-G})-(\ref{rep}) and noting that
     $$\int_{\cup_{i=1}^\ell B_{R\epsilon}(\sigma_i(x_0))}|\nabla_g\eta_\epsilon|^2dv_g=\ell\int_{B_{R\epsilon}(x_0)}|\nabla_g\eta_\epsilon|^2dv_g,$$
      we obtain
     \bea\nonumber
     \int_{\Sigma}|\nabla_g \eta_\epsilon|^2dv_g&=&\f{1}{4\pi\ell c^2}\le(
     2\log\f{1}{\epsilon}+\log(\pi\ell)-1+4\pi\ell A_{x_0}+4\pi\ell\alpha\int_\Sigma G^2dv_g\ri.\\\label{gad}
     &&\le.+O(\f{1}{R^2})+
     O(R\epsilon\log(R\epsilon))\ri).
     \eea
     Observing
     \bea\nonumber
     \int_\Sigma\eta_\epsilon dv_g&=&\f{1}{c}\le(\int_{\Sigma\setminus \cup_{i=1}^\ell B_{2R\epsilon}(\sigma_i(x_0))}Gdv_g
     +O(R\epsilon\log(R\epsilon))\ri)\\\nonumber
     &=&\f{1}{c}\le(-\int_{\cup_{i=1}^\ell B_{2R\epsilon}(\sigma_i(x_0))}Gdv_g+O(R\epsilon\log(R\epsilon))\ri)\\
     &=&\f{1}{c}O(R\epsilon\log(R\epsilon)),\label{76'}
     \eea
     we have
      \be\label{76''}\overline{\eta}_\epsilon=\f{1}{{\rm Vol}_g(\Sigma)}\int_\Sigma \eta_\epsilon dv_g=\f{1}{c}O(R\epsilon\log(R\epsilon)).\ee
      Hence
      \bna\int_\Sigma (\eta_\epsilon-\overline{\eta}_\epsilon)^2dv_g&=&\int_\Sigma \eta_\epsilon^2dv_g-2\overline{\eta}_\epsilon
      \int_\Sigma \eta_\epsilon dv_g+\overline{\eta}_\epsilon^2{\rm Vol}_g(\Sigma)\\
      &=&\f{1}{c^2}\le(\int_\Sigma G^2dv_g+O(R\epsilon\log(R\epsilon))\ri).
      \ena
       This together with (\ref{gad}) yields
       \bea\nonumber
       \|\eta_\epsilon-\overline{\eta}_\epsilon\|_{1,\alpha}^2&=&\int_\Sigma|\nabla_g\eta_\epsilon|^2dv_g-\alpha
     \int_\Sigma(\eta_\epsilon-\overline{\eta}_\epsilon)^2dv_g\\\label{norm}
       &=&\f{1}{4\pi \ell c^2}\le(
     2\log\f{1}{\epsilon}+\log(\pi\ell)-1+4\pi\ell A_{x_0}
     +O(\f{1}{R^2})+
     O(R\epsilon\log(R\epsilon))\ri).
       \eea
     Now we choose $B$ in (\ref{2pic2}) such that
     \be\label{=1}\|\eta_\epsilon-\overline{\eta}_\epsilon\|_{1,\alpha}=1.\ee
     Combining (\ref{norm}) and (\ref{=1}), we have
     \be\label{c2}
     c^2=-\f{\log\epsilon}{2\pi\ell}+\f{\log(\pi\ell)}{4\pi\ell}-\f{1}{4\pi\ell}+A_{x_0}
     +O(\f{1}{R^2})+O(R\epsilon\log(R\epsilon)).
     \ee
     It then follows from (\ref{2pic2}) and (\ref{c2}) that
     \be{\label{B}}
     B=\f{1}{4\pi\ell}+O(\f{1}{R^2})+O(R\epsilon\log(R\epsilon)).
     \ee
     Let
     \be\label{phi-e}\phi_\epsilon=\eta_\epsilon-\overline{\eta}_\epsilon.\ee
     In view of (\ref{h-g}), (\ref{phi-e}) and
     the fact that $\eta_\epsilon\in W^{1,2}(\Sigma,g)$, we have $\phi_\epsilon\in \mathscr{H}_\textbf{G}$. Moreover, the equality
     (\ref{=1}) is exactly $\|\phi_\epsilon\|_{1,\alpha}=1$, and thus (\ref{nrm}). A straightforward calculation shows on $B_{R\epsilon}(x_0)$,
     $$4\pi\ell\phi_\epsilon^2\geq 4\pi\ell
     c^2-2\log(1+\pi\ell\f{r^2}{\epsilon^2})+8\pi\ell B+O(R\epsilon\log(R\epsilon)).$$
     This together with (\ref{c2}) and (\ref{B}) yields
     $$
     \int_{B_{R\epsilon}(x_0)} e^{4\pi\ell\phi_\epsilon^2}dv_g\geq
     \pi e^{1+4\pi\ell A_{x_0}} +O(\f{1}{(\log\epsilon)^2}),
     $$
     which immediately leads to
     \be\label{BRE}\int_{\cup_{i=1}^\ell B_{R\epsilon}(\sigma_i(x_0))} e^{4\pi\ell\phi_\epsilon^2}dv_g\geq
     \pi\ell e^{1+4\pi\ell A_{x_0}} +O(\f{1}{(\log\epsilon)^2}).\ee
     Now we shall calculate the integral $\int_{\Sigma\setminus\cup_{i=1}^\ell B_{2R\epsilon}(\sigma_i(x_0))}e^{4\pi\ell\phi_\epsilon^2}dv_g$.
     By our choices of $R=-\log\epsilon$ and $c^2=O(\log\epsilon)$ (see (\ref{c2})), one can easily see that
     \be\label{c2'}R\epsilon\log(R\epsilon)=o(\f{1}{c^2}).\ee
     Recalling the representation of the Green function $G$, namely (\ref{G-representation}), one has
     \bna
     \int_{\cup_{i=1}^\ell B_{2R\epsilon}(\sigma_i(x_0))}G^2dv_g&=&\sum_{i=1}^\ell\int_{B_{2R\epsilon}(\sigma_i(x_0))}G^2dv_g\\
     &=&O((R\epsilon)^2(\log(R\epsilon))^2).
     \ena
     This together with (\ref{c2'}) gives
     \bea\nonumber
     \int_{\Sigma\setminus\cup_{i=1}^\ell B_{2R\epsilon}(\sigma_i(x_0))}G^2dv_g&=&\int_\Sigma G^2dv_g-\int_{\cup_{i=1}^\ell B_{2R\epsilon}(\sigma_i(x_0))}G^2dv_g\\\label{76''''}
     &=&\|G\|_2^2+o(\f{1}{c^2}).
     \eea
     Moreover, in view of (\ref{76'}), (\ref{76''}), (\ref{phi-e}) and (\ref{c2'}), there holds
     \bea\nonumber
     \int_{\Sigma\setminus\cup_{i=1}^\ell B_{2R\epsilon}(\sigma_i(x_0))}\phi_\epsilon^2dv_g
     &=&\int_{\Sigma\setminus\cup_{i=1}^\ell B_{2R\epsilon}(\sigma_i(x_0))}\eta_\epsilon^2dv_g+o(\f{1}{c^2})\\
     &=&\int_{\Sigma\setminus\cup_{i=1}^\ell B_{2R\epsilon}(\sigma_i(x_0))}\f{G^2}{c^2}dv_g+o(\f{1}{c^2}).\label{76'''}
     \eea
     Obviously it follows from $R=-\log\epsilon$ and (\ref{c2}) that
       \be\label{83'}\int_{\Sigma\setminus\cup_{i=1}^\ell B_{2R\epsilon}(\sigma_i(x_0))}dv_g={\rm vol}_g(\Sigma)+o(\f{1}{c^2}).\ee
     Combining (\ref{76''''})-(\ref{83'}) and using the inequality $e^t\geq 1+t$ for $t\geq 0$, we obtain
     \bea
      \label{O-BRE}\int_{\Sigma\setminus
      \cup_{i=1}^\ell B_{R\epsilon}(\sigma_i(x_0))}e^{4\pi\ell\phi_\epsilon^2}dv_g&\geq&\int_{\Sigma\setminus
      \cup_{i=1}^\ell B_{2R\epsilon}(\sigma_i(x_0))}\le(1+4\pi\ell\phi_\epsilon^2\ri)dv_g{\nonumber}\\
      &\geq& {\rm vol}_g(\Sigma)+4\pi\ell\f{\|G\|_2^2}{c^2}+o(\f{1}{c^2}).
     \eea
     Noting that $O(\f{1}{(\log\epsilon)^{2}})=o(\f{1}{c^{2}})$ and combining (\ref{BRE}) and (\ref{O-BRE}),
     we conclude (\ref{gr}) for sufficiently small $\epsilon>0$. This completes the proof of Theorem \ref{Thm1}.
     $\hfill\Box$

\section{Proof of Theorem \ref{Thm2}}\label{sec-3}
In this section, we shall prove Theorem \ref{Thm2}.
Since the proof is very similar to that of Theorem \ref{Thm1}, we only give its outline.\\

Let $j\geq 2$, $\lambda_j^\textbf{G}$ and $E_{j-1}^\perp$ be defined as in (\ref{lamj}) and (\ref{Ej}) respectively. For $\alpha<\lambda_j^\textbf{G}$, we define
\be\label{best-constant-j}\beta_j^\ast=\sup\le\{\beta:\sup_{u\in E_{j-1}^\perp,\,\|u\|_{1,\alpha}\leq 1}\int_\Sigma e^{\beta
u^2}dv_g<\infty
\ri\}.\ee
Comparing (\ref{best-constant}) with (\ref{best-constant-j}), similar to Lemma \ref{lemma-1}, we have $\beta_j^\ast=4\pi\ell$, where $\ell$ is defined as in (\ref{ell}).

We now prove $(ii)$ of Theorem \ref{Thm2}.
If $\alpha\geq \lambda_j^\textbf{G}$ and $\beta>0$, we take $u_j\in\mathscr{H}_\textbf{G}\cap C^1(\Sigma,g)$
satisfies $\Delta_gu_j=\lambda_{j}^\textbf{G} u_j$ and $u_j\not\equiv 0$.
It follows that \be\label{en-2}\int_\Sigma|\nabla_g(tu_j)|^2dv_g-\alpha\int_\Sigma (tu_j)^2dv_g\leq 0,\quad \forall t\in\mathbb{R}\ee
and that
\be\label{tend}\int_\Sigma e^{\beta(tu_j)^2}dv_g\ra\infty\quad {\rm as}\quad t\ra \infty.\ee
Then (\ref{en-2}) and (\ref{tend}) imply that the supremum in (\ref{ineq-2}) is infinity.

If $\alpha<\lambda_j^\textbf{G}$ and $\beta>4\pi\ell$, then we shall prove that the supremum in (\ref{ineq-2}) is infinity.
To do this, we let $\{e_i\}_{i=1}^{m_{j-1}}\subset \mathscr{H}_\textbf{G}\cap C^1(\Sigma,g)$
be an orthonormal basis of $E_{j-1}=E_{\lambda_1^\textbf{G}}\oplus\cdots\oplus E_{\lambda_{j-1}^\textbf{G}}$
with respect to the inner product on $L^2(\Sigma,g)$,
namely, $E_{j-1}={\rm span}\{e_1,\cdots,e_{m_{j-1}}\}$ and
$$(e_i,e_k)=\int_\Sigma e_ie_kdv_g=\delta_{ik}=\le\{\begin{array}{lll}
1,&i=k\\[1.5ex]
0,&i\not=k
\end{array}\ri.$$
for all $i,k=1,\cdots,m_{j-1}$.
 Let $\widetilde{M}_k$
be defined as in (\ref{test}). Set
$$Q_k=\widetilde{M}_k-\f{1}{{\rm Vol}_g(\Sigma)}\int_\Sigma \widetilde{M}_kdv_g-\sum_{i=1}^{m_{j-1}}(\widetilde{M}_k,e_i) e_i.$$
 Then $Q_k\in E_{j-1}^\perp$. By a straightforward calculation,
$\|Q_k\|_{1,\alpha}^2=(1+O(r))8\pi\ell\log k+O(1)$. Denote $Q_k^\ast=Q_k/\|Q_k\|_{1,\alpha}$.
It follows that for any fixed $\beta>4\pi\ell$,
\bna
\int_{\Sigma}e^{\beta {Q_k^\ast}^2}dv_g&\geq&\int_{B_{rk^{-1/4}}(x_0)}e^{
\f{\beta(\log k+O(1))^2}{(1+O(r))8\pi\ell\log k+O(1)}}dv_g\\
&=&e^{\f{\beta(1+o_k(1))\log k}{(1+O(r))8\pi\ell}}\pi r^2k^{-1/2}(1+o_k(1)).
\ena
Choosing $r>0$ sufficiently small and then passing to the limit $k\ra\infty$ in the above estimate, we
conclude
$$\int_{\Sigma}e^{\beta {Q_k^\ast}^2}dv_g\ra \infty\quad{\rm as}\quad k\ra\infty.$$
Hence the supremum in (\ref{ineq-2}) is infinity. \\

In the following, we sketch the proof of $(i)$ and $(iii)$ of Theorem \ref{Thm2}.

 Let $\alpha<\lambda_j^\textbf{G}$. By a direct method of variation, one can see that for any $0<\epsilon<4\pi\ell$,
there exists some $u_\epsilon\in E_{j-1}^\perp$ with $\|u_\epsilon\|_{1,\alpha}=1$ such that
$$\int_\Sigma e^{(4\pi\ell-\epsilon)u_\epsilon^2}dv_g=\sup_{u\in E_{j-1}^\perp,\|u\|_{1,\alpha}\leq 1}
\int_\Sigma e^{(4\pi\ell-\epsilon)u^2}dv_g.$$
Clearly $u_\epsilon$ satisfies the Euler-Lagrange equation
 \be\label{E-L-2}\le\{\begin{array}{lll}
\Delta_gu_\epsilon-\alpha u_\epsilon=\f{1}{\lambda_\epsilon}u_\epsilon e^{(4\pi\ell-\epsilon)u_\epsilon^2}
-\f{\mu_\epsilon}{\lambda_\epsilon}-\sum_{k=1}^{m_{j-1}}\gamma_ke_k\\[1.5ex]
u_\epsilon\in E_{j-1}^\perp,\,\,\|u_\epsilon\|_{1,\alpha}=1\\[1.5ex]
\lambda_\epsilon=\int_\Sigma u_\epsilon^2e^{(4\pi\ell-\epsilon)u_\epsilon^2}dv_g\\[1.5ex]
\mu_\epsilon=\f{1}{{\rm Vol}_g(\Sigma)}\int_\Sigma u_\epsilon e^{(4\pi\ell-\epsilon)u_\epsilon^2}dv_g\\[1.5ex]
\gamma_k=\int_\Sigma\f{1}{\lambda_\epsilon}e_ku_\epsilon e^{(4\pi\ell-\epsilon)u_\epsilon^2}dv_g.
\end{array}\ri.\ee
Without loss of generality, we assume $c_\epsilon=u_\epsilon(x_\epsilon)=\sup_\Sigma |u_\epsilon|\ra +\infty$ and $x_\epsilon\ra x_0$
as $\epsilon\ra 0$. Let $r_\epsilon$ be the blow-up scale defined as in (\ref{scale}) and
$\varphi_\epsilon(y)=c_\epsilon(u_\epsilon(\exp_{x_\epsilon}(r_\epsilon y))-c_\epsilon)$ for $y\in \mathbb{B}_{\delta r_\epsilon^{-1}}(0)$,
where $0<\delta<\f{1}{2}i_g(\Sigma)$.  As before, we can derive
$$\varphi_\epsilon(y)\ra \varphi(y)=-\f{1}{4\pi\ell}
\log(1+\pi\ell|y|^2)\quad{\rm in}\quad C^1_{\rm loc}(\mathbb{R}^2).$$
Moreover, we can prove that  $\forall 1<q<2$, $c_\epsilon u_\epsilon$ converges to a Green function¡¡
$G$ weakly in $W^{1,q}(\Sigma,g)$, strongly in $L^{\f{2q}{2-q}}(\Sigma,g)$, and almost everywhere in $\Sigma$.
In this case, $G$ satisfies
\be\label{Green-2}\le\{\begin{array}{lll}
 \Delta_gG-\alpha G=\f{1}{\ell}\sum_{i=1}^\ell \delta_{\sigma_i(x_0)}-\f{1}{{\rm Vol}_g(\Sigma)}-\sum_{j=1}^{m_{j-1}}e_k(x_0)e_k\\[1.5ex]
 \int_\Sigma G\phi dv_g=0,\,\,\,\forall \phi\in E_{j-1}\\[1.5ex]
 G(\sigma(x))=G(x),\,\forall x\in \Sigma\setminus \textbf{G}(x_0),\,\forall \sigma\in \textbf{G}.
 \end{array}\ri.\ee
 Similarly, $G$ has a decomposition (\ref{G-representation}) near $x_0$ and $A_{x_0}$ is defined as in (\ref{ax00}).
 Analogous to Proposition \ref{up-bd}, we arrive at
 \be\label{u-bb}\sup_{u\in E_{j-1}^\perp,\,\|u\|_{1,\alpha}\leq 1}\int_\Sigma e^{4\pi\ell u^2}dv_g
\leq {\rm Vol}_g(\Sigma)+
\pi\ell e^{1+4\pi\ell A_{x_0}}.\ee
This particularly leads to $(i)$ of Theorem \ref{Thm2}.

Finally we construct a sequence of functions to show that the estimate (\ref{u-bb}) is not true. This implies that blow-up can not occur and elliptic
estimates on (\ref{E-L-2}) give a desired extremal function.
 To do this, we let $\eta_\epsilon$, $\phi_\epsilon$ be defined respectively as in (\ref{eta-e}) and (\ref{phi-e}) satisfying
 $\eta_\epsilon\in W^{1,2}(\Sigma,g)$ and $\|\phi_\epsilon\|_{1,\alpha}=1$. Note that the constants $c$ and $B$ in definitions of
 $\eta_\epsilon$ and $\phi_\epsilon$ are given by (\ref{c2}) and (\ref{B}) respectively. It then follows that
 \be\label{gr-3}\int_{\Sigma}e^{4\pi\ell \phi_\epsilon^2}dv_g\geq {\rm Vol}_g(\Sigma)+\pi\ell e^{1+4\pi\ell A_{x_0}}+
 \f{4\pi\ell \|G\|_{L^2(\Sigma,g)}^2}{-\log\epsilon}+o(\f{1}{-\log\epsilon}).\ee
  Let
 \be\label{t-2}\widetilde{\phi}_\epsilon=\phi_\epsilon-\sum_{i=1}^{m_{j-1}}(\phi_\epsilon,e_i)e_i.\ee
 Obviously $\widetilde{\phi}_\epsilon\in E_{j-1}^\perp$. Since $G$ satisfies (\ref{Green-2}) and
 $$\int_\Sigma Ge_idv_g=\lim_{\epsilon\ra 0}\int_{\Sigma}c_\epsilon u_\epsilon e_i dv_g=0,\quad\forall 1\leq i\leq m_{j-1},$$
 we calculate
 \bna
 (\phi_\epsilon,e_i)=\int_{\cup_{i=1}^\ell B_{2R\epsilon}(\sigma_i(x_0))} (\eta_\epsilon-\overline{\eta}_\epsilon)e_idv_g
 +\int_{\Sigma\setminus\cup_{i=1}^\ell B_{2R\epsilon}(\sigma_i(x_0))}\le(\f{G}{c}-\overline{\eta}_\epsilon\ri)e_idv_g=o(\f{1}{\log^2\epsilon}).
 \ena
 This together with (\ref{t-2}) leads to
 \be\label{p-t}\widetilde{\phi}_\epsilon=\phi_\epsilon+o(\f{1}{\log^2\epsilon}),\quad
 \|\widetilde{\phi}_\epsilon\|_{1,\alpha}^2=1+o(\f{1}{\log^2\epsilon}).\ee
 It follows from (\ref{gr-3}) and (\ref{p-t}) that
 \bna
 \int_\Sigma e^{4\pi\ell\f{\widetilde{\phi}_\epsilon^2}{\|\widetilde{\phi}_\epsilon\|_{1,\alpha}^2}}dv_g&=&
 \int_\Sigma e^{4\pi\ell \phi_\epsilon^2+o(\f{1}{-\log\epsilon})}dv_g\\
 &\geq&\le(1+o(\f{1}{-\log\epsilon})\ri)\le({\rm Vol}_g(\Sigma)+\pi\ell e^{1+4\pi\ell A_{x_0}}+
 \f{4\pi\ell \|G\|_{L^2(\Sigma,g)}^2}{c^2}+o(\f{1}{c^2})\ri)\\
 &\geq& {\rm Vol}_g(\Sigma)+\pi\ell e^{1+4\pi\ell A_{x_0}}+
 \f{4\pi\ell \|G\|_{L^2(\Sigma,g)}^2}{-\log\epsilon}+o(\f{1}{-\log\epsilon}),
 \ena
 which implies that (\ref{u-bb}) does not hold. This completes the proof of $(iii)$ of Theorem \ref{Thm2}.

\bigskip

{\bf Acknowledgements}. The authors thank the referee for very careful reading and constructive comments, which greatly improve
 the previous version of this paper. This work is partly supported by the National Science Foundation of China (Grant Nos. 11171347 and
  11471014).


\begin{thebibliography}{00}

\bibitem{A-D} A. Adimurthi, O. Druet, Blow-up analysis in dimension 2 and a sharp form of Moser-Trudinger inequality, Comm. Partial Differential Equations 29 (2004) 295-322.

\bibitem{A-S} A. Adimurthi, M. Struwe, Global compactness
properties of semilinear elliptic equation with critical exponential
growth, J. Functional Analysis 175 (2000) 125-167.


\bibitem{C-C} L. Carleson, A. Chang, On the existence
of an extremal function for an inequality of J. Moser, Bull. Sci.
Math. 110 (1986) 113-127.

\bibitem{Chen-90} W. Chen, A Trudinger inequality on surfaces with conical singularities,
Proc. Amer. Math. Soc. 108 (1990) 821-832.

\bibitem{C-L} W. Chen and C. Li, Classification of solutions of some
nonlinear elliptic equations, Duke Math. J.  63 (1991) 615-622.

\bibitem{de-doO} M. de Souza, J. M. do \'O, A sharp Trudinger-Moser type inequality in $\mathbb{R}^2$, Trans. Amer. Math. Soc.
366 (2014) 4513-4549.


\bibitem{DJLW} W. Ding, J. Jost, J. Li, G.Wang, The differential equation $-\D u=8\pi-8\pi he^u$ on a compact Riemann Surface, Asian
J. Math. 1 (1997) 230-248.


\bibitem{Flucher} M. Flucher,  Extremal functions for Trudinger-Moser
inequality in 2 dimensions, Comment. Math. Helv. 67 (1992)
471-497.

\bibitem{Fontana} L. Fontana, Sharp bordline Sobolev inequalities on compact Riemannian manifolds, Comment. Math. Helv. 68
(1993) 415-454.


\bibitem{Iula-Mancini} S. Iula, G. Mancini, Extremal functions for singular Moser-Trudinger embeddings, Nonlinear Anal. 156 (2017) 215-248.

\bibitem{Lijpde} Y. Li, Moser-Trudinger inequality on compact
Riemannian manifolds of dimension two, J. Partial Differential
Equations 14 (2001) 163-192.

\bibitem{Liscience} Y. Li, The existence of the extremal function
of Moser-Trudinger inequality on compact Riemannian manifolds,
Science in China (Series A) 48 (2005) 618-648.


\bibitem{Lin} K. Lin, Extremal functions for Moser's
inequality, Trans. Amer. Math. Soc. 348 (1996) 2663-2671.

\bibitem{Lu-Yang-1} G. Lu, Y. Yang, The sharp constant and extremal functions for Moser-Trudinger inequalities involving $L^{p}$ norms, Discrete and Continuous Dynamical Systems 25 (2009) 963-979.

\bibitem{Lu-Yang-2} G. Lu, Y. Yang, Adams' inequalities for bi-Laplacian and extremal functions in dimension four, Adv. Math. 220 (2009) 1135-1170.

\bibitem{Mancini-Martinazzi} G. Mancini, L. Martinazzi, The Moser-Trudinger inequality and its extremals on a disk via energy estimates, Calc. Var. (2017)
56:94,  https://doi.org/10.1007/s00526-017-1184-y.

\bibitem{Moser}  J. Moser, A sharp form of an inequality by N. Trudinger, Indiana Univ. Math. J. 20 (1971) 1077-1091.

\bibitem{Moser-73} J. Moser, On a nonlinear problem in differential geometry, Dynamical systems (Proc. Sympos., Univ. Bahia, Salvador, 1971), 273-280, Academic Press, New York, 1973.

\bibitem{Nguyen} V. H. Nguyen, Extremal functions for the Moser-Trudinger inequality of Adimurthi-Druet type in $W^{1,N}(\mathbb R^N)$, arXiv: 1702.07970.

\bibitem{Peetre} J. Peetre, Espaces d'interpolation et theoreme de Soboleff, Ann. Inst. Fourier (Grenoble) 16 (1966) 279-317.

\bibitem{Pohozaev} S. Pohozaev, The Sobolev embedding in the special case
$pl=n$, Proceedings of the technical scientific conference on
advances of scientific reseach 1964-1965, Mathematics sections,
158-170, Moscov. Energet. Inst., Moscow, 1965.


\bibitem{Struwe} M. Struwe, Critical points of embeddings of $H_0^{1,n}$ into Orlicz spaces, Ann. Inst. H. Poincar\'e,
Analyse Non Lin\'eaire 5 (1988) 425-464.


\bibitem{Tintarev} C. Tintarev, Trudinger-Moser inequality with remainder terms, J. Funct. Anal. 266 (2014) 55-66.


\bibitem{Trudinger} N. Trudinger, On embeddings into Orlicz spaces and
some applications, J. Math. Mech. 17 (1967) 473-484.



\bibitem{InterJM} Y. Yang, Extremal functions for a sharp Moser-Trudinger inequality, Inter. J. Math. 17 (2006) 331-338.

\bibitem{Yang-JFA-06} Y. Yang, A sharp form of Moser-Trudinger inequality in high dimension, J. Funct. Anal. 239 (2006) 100-126.


\bibitem{Yang-Tran} Y. Yang, A sharp form of the Moser-Trudinger inequality on a compact Riemannian surface,
Trans. Amer. Math. Soc. 359 (2007) 5761-5776.


\bibitem{Yang-JDE-15} Y. Yang, Extremal functions for Trudinger-Moser inequalities of Adimurthi-Druet type in dimension two,
J. Differential Equations 258 (2015) 3161-3193.

\bibitem{YangJGA} Y. Yang, A Trudinger-Moser inequality on compact Riemannian surface
  involving Gaussian curvature, J. Geom. Anal. 26 (2016) 2893-2913.
\bibitem{Yang-Zhu-JFA} Y. Yang, X. Zhu, Blow-up analysis concerning singular Trudinger-Moser inequalities in dimension two, J. Funct. Anal. 272 (2017) 3347-3374.

\bibitem{Yang-Zhu-SCM} Y. Yang, X. Zhu, Existence of solutions to a class of Kazdan-Warner equations on compact Riemannian surface,
Sci. China Math. 61 (2018) 1109-1128.

\bibitem{YuanHuang} A. Yuan, Z. Huang, An improved singular Trudinger-Moser inequality in dimension two,
Turkish J. Math. 40 (2016) 874-883.

\bibitem{YuanZhu} A. Yuan, X. Zhu, An improved singular Trudinger-Moser inequality in unit ball, J. Math. Anal. Appl. 435
(2016) 244-252.


\bibitem{Yudovich} V. Yudovich, Some estimates connected with integral operators and with solutions of elliptic equations,
Sov. Math. Docl. 2 (1961) 746-749.

\bibitem{ZhuJ} J. Zhu, Improved Moser-Trudinger inequality involving $L^p$ norm in $n$ dimensions, Advanced Nonlinear Study
14 (2014) 273-293.

\end{thebibliography}
\end{document}